\documentclass[12pt,reqno]{amsart}
\usepackage{amssymb,delarray}
\usepackage{amsfonts,amscd}
\usepackage{epsfig}
\usepackage[all]{xy}

\def\leq{\leqslant}
\def\geq{\geqslant}

\def\pr{{\sf pr}}

\def\pr{\mathsf{pr}}
\newcommand{\Sing}{\operatorname{Sing}}

\newtheorem{thm}{Theorem}
\newtheorem{lem}
{Lemma}
\newtheorem{prop}
{Proposition}
\newtheorem{prob}
{Problem}
{Claim}
\newtheorem{df}
{Definition}
\newtheorem{cor}
{Corollary}
\newtheorem{rem}
{Remark}
{Question}
\newtheorem{ex-thm}{Theorem-Example}

{\catcode`\@=11
\gdef\n@te#1#2{\leavevmode\vadjust{%
 {\setbox\z@\hbox to\z@{\strut#1}%
  \setbox\z@\hbox{\raise\dp\strutbox\box\z@}\ht\z@=\z@\dp\z@=\z@%
  #2\box\z@}}}
\gdef\leftnote#1{\n@te{\hss#1\quad}{}}
\gdef\rightnote#1{\n@te{\quad\kern-\leftskip#1\hss}{\moveright\hsize}}
\gdef\?{\FN@\qumark}
\gdef\qumark{\ifx\next"\DN@"##1"{\leftnote{\rm##1}}\else
 \DN@{\leftnote{\rm??}}\fi{\rm??}\next@}}

\begin{document}
\baselineskip=13.7pt plus 2pt 

\title[On germs of finite morphisms of smooth surfaces]{On germs of finite morphisms of smooth surfaces}

\author[Vik.S. Kulikov]{Vik.S. Kulikov}

\address{Steklov Mathematical Institute of Russian Academy of Sciences, Moscow, Russia}
 \email{kulikov@mi.ras.ru}

\dedicatory{} \subjclass{}


\maketitle

\def\st{{\sf st}}

\quad \qquad \qquad

\begin{abstract} Questions related to deformations of germs of finite morphisms of smooth surfaces are discussed.
A classification of the  four-sheeted germs of finite covers $F: (U,o')\to (V,o)$ is given up to smooth deformations, where  $(U,o')$ and $(V,o)$ are two connected germs of smooth complex-analytic surfaces. The singularity types of their branch curves and the local monodromy groups are investigated also.
\end{abstract} \vspace{0.5cm}

\def\st{{\sf st}}


\setcounter{section}{-1}

\section{Introduction}
Let $(V,o)=(\mathbb B_{\varepsilon},o)$ be a ball in $\mathbb C^2$ of small radius $\varepsilon>0$, $(U,o')$ a connected germ of smooth complex-analytic surface, and $F: (U,o')\to (V,o)$ a germ of finite holomorphic mapping (below,  a {\it germ of finite cover}) of local degree $\deg_{o'} F=d$, given in local coordinates $z,w$ in $(U,o')$ and $u,v$ in $(V,o)$ by functions
$$  u= f_1(z,w), \qquad  v= f_2(z,w), $$
where $f_i(z,w)\in \mathbb C[[z,w]]$ are convergent power series. Denote by $R\subset (U,o')$ the ramification divisor of $F$ given by equation
$$J(f):= \det \left(\begin{array}{cc} \frac{\partial u}{\partial z} & \frac{\partial u}{\partial w} \\ \frac{\partial v}{\partial z} & \frac{\partial v}{\partial w}\end{array}\right)= 0 $$
and by $B=F(R_{red})\subset (V,o)$ the germ of branch curve of $F$.

The germ $F$ defines a homomorphism $F_*:\pi_1(V\setminus B,p)\to \mathbb S_d$ (the {\it monodromy} of the germ $F$), where $\mathbb S_d$ is the symmetric group acting on the fibre $F^{-1}(p)$.  The group $G_F=\text{im} F_*$ is called the ({\it local}) {\it monodromy group} of $F$. Note that $G_F$ is a transitive subgroup of $\mathbb S_d$.

We say that two germs $F_1: (U,o')\to (V,o)$ and $F_2: (U,o')\to (V,o)$ are {\it equivalent} if they differ from each other by changes of coordinates in $(U,o')$ and $(V,o)$.

The aim of this article is to investigate the germs of finite covers with up to {\it deformation equivalence}. In short (see Definition \ref{def-d} in  Subsection \ref{defm}), two germs $F_1$ and $F_2$ are deformation equivalent if they can be included in a smooth family of the germs of  finite covers preserving the singularity type of germs of branch curves.

The germs of finite covers of local degree $\deg_{o'} F=3$ were investigated in \cite{Ku-3}. In this article, we investigate the singularity types of the germs of  branch curves  of  finite covers of local degree $\deg_{o'} F=4$, describe their monodromy groups, and give a complete classification of them with up to deformation equivalence.

Note that, a priori, the monodromy group $G_F$ of the four-sheeted germ of a finite cover $F$ is one of the following subgroups of $\mathbb S_4$:
a cyclic group $\mathbb Z_4$ generated by a cycle of length 4, the Klein four group $Kl_4$, a dihedral group $\mathbb D_4$, the alternating group $\mathbb A_4$, and $\mathbb S_4$.

Before to formulate the main result of the article, we introduce several notations.
\begin{itemize}
\item $T[f(u,v)=0]=T(B)$, the singularity type of a curve germ  $B$ given by equation $f(u,v)=0$;
\item $A_n:=T[v^2-u^{n+1}=0]$, $n\geq 0$;
\item $D_n:=T[u(v^2-u^{n-2})=0]$, $n\geq 4$;
\item $E_6:=T[v(v^2-u^3)=0]$;
\item $E_7:=T[v^3-u^4=0]$;
\item $E_8:=T[v^3-u^5=0]$;
\item $T_3(n,\beta):=T[v^3-u^{3n+\beta}=0]$, $n\geq 0$, $\beta=1,2$;
\item $T_{3,2k+n,2k}:=T[v((v-u^k)^2-u^{2k+n+1})=0]$, $k\geq 1$, $n\geq 0$;
\item $T_{3,n-1,n}:= T[v(v^2-u^{n})=0]$, $n\geq 2$;
\item $T_{4, 2n_1,2n_2}:=[(v^2-u^{2n_1+1})(u^2-v^{2n_2+1})=0]$, $n_1,n_2\in\mathbb N$.
\end{itemize}

Note that in these notations there are several intersections, for example, $E_6=T_{3,2,3}$ and $E_7=T_3(1,1)$, $E_8=T_3(1,2)$.

The main result of the article is the following theorem in which we give a complete classification of four-sheeted germs of finite  covers with up to deformation equivalence and compute their main invariants: the singularity types of branch curves and their monodromy groups.
\begin{thm} \label{main} A germ of finite cover $F:(U,o')\to (V,o)$, $\deg_{o'} F=4$, given by functions $u=f_1(z,w)$ and $v=f_2(z,w)$, either is equivalent to one of the following germs of finite covers:
\begin{itemize}
\item[$({\bf 1}_1)$]  $F_{{\bf 1}_1}= \{ u=z,\quad v=w^4\} $, $G_{F_{{\bf 1}_1}}=\mathbb Z_4$, $T(B_{{\bf 1}_1})=A_0$;
\item[$({\bf 2}_1)$]  $F_{{\bf 2}_1}=\{ u=z^2,\quad v=w^2\}$,  $G_{F_{{\bf 2}_1}}=Kl_4$,\,\,\, $T(B_{{\bf 2}_1})=A_1$;
\item[$({\bf 3}_1)$]  $F_{{\bf 3}_1,n}= \{  u=z,\quad v= w^4-2z^nw^2\} $, $G_{F_{{\bf 3}_1,n}}=\mathbb D_4$,
$T(B_{{\bf 3}_1,n})=A_{4n-1}$, $n\geq 1$;
\item[$({\bf 3}_2)$]  $F_{{\bf 3}_2,n}= \{  u=z^2,\quad v= w^2-z^{2n+1}\} $, $G_{F_{{\bf 3}_2,n}}=\mathbb D_4$,
$T(B_{{\bf 3}_2,n})=D_{2n+3}$, $n\geq 1$;
\item[$({\bf 4}_1)$]  $F_{{\bf 4}_1,n}= \{ u=z, \quad v=w^4+w^3z^{n}\} $,
$G_{F_{{\bf 4}_1,n}}=\mathbb S_4$, $T(B_{{\bf 4}_1,n})=A_{8n-1}$, $n\geq 1$,
\end{itemize}
or is deformation equivalent  to one of the following germs of finite covers:
\begin{itemize}
\item[$({\bf 4}_2)$]  $F_{{\bf 4}_2,n,\beta}= \{ u=z, \qquad v=w^4+4wz^{3n+\beta}\} $,\,\,
$G_{F_{{\bf 4}_2,n,\beta}}=\mathbb S_4$,\,\, $T(B_{{\bf 4}_2,n,\beta})=\newline T_3(4n+\beta,\beta)
$, $n\geq 0$, $\beta=1,2$; 
\item[$({\bf 4}_3)$] $F_{{\bf 4}_3,n,m}= \{ u=z,\,\, v = w^4-\frac{8}{3}z^mw^3+2(z^{2m}-z^{2(m+n)+1})w^2\} $, $G_{F_{{\bf 4}_3,n,m}}=\mathbb S_4$, $T(B_{F_{{\bf 4}_3,n,m}})= T_{3,8m+6n+2, 8m}$, $m\geq 1$, $n\geq 0$;
\item[$({\bf 4}_4)$] $F_{{\bf 4}_4,n,m}= \{ u=z,\,\, v = w^4 +12z^{m+n+1}w^3 -2z^{2m+1}w^2\}$,  $G_{F_{{\bf 4}_4,n,m}}=\mathbb S_4$, $T(B_{F_{{\bf 4}_4,n,m}})= T_{3,8m+2n+4, 8m+4}$, $m\geq 1$, $n\geq 0$;
\item[$({\bf 4}_5)$] $F_{{\bf 4}_5,k,m}= \{ u=z,\,\, v = w^4-\frac{8}{3}z^m(1+2z^k)w^3+2z^{2m}(1+4z^k+3z^{2k})w^2\}$,
$G_{F_{{\bf 4}_5,k,m}}=\mathbb S_4$, $T(B_{F_{{\bf 4}_5,k,m}})= T_{3,8m+6k-1, 8m}$, $m\geq 1$, $k\geq 1$;
\item[$({\bf 4}_6)$] $F_{{\bf 4}_6,n,m}= \{ u=z,\,\, v = w^4-\frac{4}{3}z^m(1+3z^n)w^3+4z^{2m}(1+z^n)w^2\} $,
$G_{F_{{\bf 4}_6,n,m}}=\mathbb S_4$, $T(B_{F_{{\bf 4}_6,n,m}})= T_{3,8m+2n-1, 8m}$, $m\geq 1$, $n\geq 1$;
\item[$({\bf 4}_7)$] $F_{{\bf 4}_7,m}= \{ u=z,\,\, v = w^4-\frac{16}{3}z^mw^3+6z^{2m}w^2 \} $, $G_{F_{{\bf 4}_7,m}}=\mathbb S_4$, $T(B_{F_{{\bf 4}_7,m}})= T_{3,8m-1, 8m}$, $m\geq 1$;
\item[$({\bf 4}_8)$] $F_{{\bf 4}_8,n_1,n_2}= \{ u=  z^2+w^{2n_1+1}, \quad v =  w^2+z^{2n_2+1}\} $, $G_{F_{{\bf 4}_8n_1,n_2}}=\mathbb S_4$, $T(B_{F_{{\bf 4}_8,n_1,n_2}})= T_{4,2n_1, 2n_2}$, $n_1,n_2\geq 1$;
\end{itemize}

In case $({\bf 4}_2)$, any two germs of finite covers deformation equivalent to the germ $F_{{\bf 4}_2,0,1}$ are equivalent.
\end{thm}

\begin{cor} \label{cor1} The set of monodromy groups $G_F$ of the four-sheeted germs of finite covers is
$\{ \mathbb Z_4, Kl_4, \mathbb D_4, \mathbb S_4\}.$

The alternating group $\mathbb A_4$ is not the monodromy group of any four-sheeted germ of finite cover.
\end{cor}

Moreover, in Proposition \ref{alt4} it is proved that the alternating group $\mathbb A_4$ is not the monodromy group of any germ of  finite cover, nevertheless (see Proposition \ref{alt2n-1}), the groups $\mathbb A_{2n-1}$ for $n\geq 2$ are the monodromy groups of germs of finite covers.

Note that all possible deformation types of four-sheeted germs of finite covers, described in Theorem \ref{main}, with the exception of two (germs $F_{{\bf 4}_5,k,m}$ and $F_{{\bf 4}_6,3k,m}$), have different pair of main invariants.
The monodromy groups of the germs $F_{{\bf 4}_5,k,m}$ and $F_{{\bf 4}_6,3k,m}$ are the same and the singularity types of their branch curves are the same also, but deformation types of these covers differ in the singularity types of their ramification divisors.
Therefore, the following problem is interesting.

\begin{prob} To find a finite set of invariants defining completely the deformation types of the germs of finite covers.\end{prob}
In particular, would be interesting to know if there are deformation non-equivalent germs of finite covers of local degree $d\geq 5$, the singularity types of branch curves and ramification divisors of which are the same, the covers have conjugated monodromy groups, and the sets of conjugacy classes of the images of geometric generators of the local fundamental groups of the branch curves are the same also.

In Section 1, the main definitions are given and several auxiliary lemmas are proved. In this section, we consider different aspects of the relationship between the deformation equivalence of germs of finite covers, the singularity types of their branch curves, and the monodromy groups of these germs.  The proof of Theorem \ref{main} is given in Section 2. In Section 3, we give examples of finite groups that are not the monodromy groups of any germs of finite covers.

\section{Definitions and preliminary facts}
\subsection{Equisingular deformations of curve germs.}\label{d-loc}
Let $$\mathbb B_{r}=\{ (u,v)\in\mathbb C^2 \mid |u|^2+|v|^2<r^2 \}$$ be a ball of radius $r$ with center at the point $o$ and $C\subset \mathbb B_{r}$ a reduced 
curve given by equation $h(u,v)=0$, where $h(u,v)\in H^0(\mathbb B_{r},
\mathcal O_{\mathbb B_r})$. Below, the words {\it "$B$ is a curve germ in $(V,o)$"}\, means that $V=\mathbb B_{\varepsilon}$, were $\varepsilon\ll r$, and $B=C\cap V$.

Let $B_1,\dots , B_k$ be the irreducible components of a curve germ $B\subset V$ and $\sigma=\sigma_1\circ \dots \circ \sigma_n: V_n\to V$ the minimal sequence of $\sigma$-processes $\sigma_i:V_i\to V_{i-1}$ with centers at points resolving the singular point $o$ of $B$ and such that $\sigma^{-1}(B)$ is a divisor with normal crossings. Denote by $E_{i+k}\subset V_n$ the proper inverse image of the exceptional curve of $\sigma_i$ and $B_j'\subset V_n$ the proper inverse image of $B_j$.
Put $\delta_{i,j}:=(E_i,E_j)_{V_n}$ for $k+1\leq i,j\leq k+n$, $i\neq j$; $\delta_{i,j}:=(B'_i,E_j)_{V_n}$ for $1\leq i\leq k$ and $k+1\leq j\leq k+n$; $\delta_{i,j}:=0$ for $1\leq i,j\leq k$, and  $\delta_{i,i}:=0$ for $1\leq i\leq k+n$. Note that, by definition,  $\delta_{i,j}$ takes only two values: $0$ or $1$.

Let $\Gamma(B)$ be the graph of the curve germ $B$.

\begin{df} \label{def-gr} The graph $\Gamma(B)$ is a weighted graph having $n+m$ vertices $v_i$. The vertices $v_i:=b_i$, $i=1,\dots, k$,  are in one-to-one correspondence with the curve germs $B'_1,\dots, B'_m$ and they have weights $w_i=0$. The vertices $v_{i+k}:=e_{i+m}$, $i=1,\dots , n$,  are in one-to-one correspondence with the curves $E_{1+k},\dots, E_{n+k}$ and they have weights $w_{i+k}=(E_{i+k}^2)_{V_n}$.  The vertices $v_i$ and $v_j$ are connected by the edge of $\Gamma(B)$ if and only if $\delta_{i,j}=1$.
\end{df}

Note that the graph $\Gamma(B)$ is a tree if $(B,o)$ is not a divisor with normal crossings consisting of two components.

By definition, a {\it family of curve germs}  is a triple $(\mathcal V=V\times D_{\varepsilon},
\mathcal B, \text{pr}_2)$, where $V=\mathbb B_{r}\subset \mathbb C^2$ is a ball,
$D_{\varepsilon}=\{ t\in \mathbb C\mid |t|<\varepsilon\}$ is a disc in the complex plane, $\mathcal B$ is an effective reduced
divisor in $\mathcal V$, the  restriction to $\mathcal B$ of the projection $\text{pr}_2:\mathcal V\to D_{\varepsilon}$
 is a flat holomorphic map.

\begin{df} \label{def1} {\rm (\cite{Z1}, 
see also \cite{W})}
A family $(\mathcal V, \mathcal B, \text{pr}_2)$ is an {\it equisingular deformation} of curve germs if
$\Sing\, \mathcal B=\{ o\}\times D_{\varepsilon}$
and there exists a finite sequence of  $n$ monoidal transformations {\rm (}blowups{\rm )}  $\widetilde{\sigma}_i:\mathcal V_{i}\to \mathcal V_{i-1}$ {\rm (}where $\mathcal V_0=\mathcal V${\rm )}  with centers in smooth curves $\mathcal S_i\subset \Sing \,  \mathcal B_i$, where  $\mathcal B_0=\mathcal B$ and $\mathcal B_{i+1}=\widetilde{\sigma}_{i}^{-1}(\mathcal B_i)$, and such that
\begin{itemize}
\item[$(i)$] $Sing\, \mathcal B_i$ is a disjoint union of sections of $\text{pr}_2\circ\widetilde{\sigma}_1\circ\dots\circ\widetilde{\sigma}_{i-1}$  for each $i$,
\item[$(ii)$] $\mathcal B_{n}$ is a divisor with normal crossings in  $\mathcal V_{n}$.
\end{itemize}
\end{df}
Note that in this case the divisor $\mathcal B_{n}$ has only simple double points as its singular points.

We say that two curve germs are {\it strongly equisingular equivalent} if they can be imbedded in an equisingular deformation of curve germs as fibres of $\text{pr}_2$. Let us continue the  strong equisingular equivalence to an equivalence relation and will say that two curve germs have the same {\it singularity type} iff they are equisingular equivalent.


The following two propositions are well known.
\begin{prop} {\rm (}\cite{W}{\rm )} Two curve germs $(B_1,o)$ and $(B_2,o)$ are equisingular equivalent if and only if their graphs $\Gamma(B_1)$ and $\Gamma(B_2)$ are isomorphic as weighted graphs.
\end{prop}

We  say that a curve germ $(B,o)\subset (V,o)$ is {\it rigid} if for any curve germ $(B_1,o)\subset (V,o)$ equisingular equivalent to $(B,o)$, there is a biholomorphic mapping $G: (V,o)\to (V,o)$ such that $G(B_1)=B$.

\begin{prop} {\rm (\cite{Ar})} \label{Anrig} Curve germs $(B,o)$ having one of the following singularity types $A_n$, $D_n$, $E_6,$, $E_7$, $E_8$ are rigid.
\end{prop}

The  imbedding  $\mathbb B_{r_2}\subset \mathbb B_{r_1}\subset \mathbb B_{r}$ of balls, $r_2\leq r_1\leq r$, induces a homo\-morphism
$i_*:\pi_1(\mathbb B_{r_2}\setminus B)\to\pi_1(\mathbb B_{r_1}\setminus B)$ of the fundamental groups.
The following theorem is well known (see, for example, \cite{HLe}).

\begin{thm} \label{loc-pi}
There is  a radius $r_0$ such that for any $\varepsilon\leq r_0$ the homomorphism
$i_*:\pi_1(\mathbb B_{\varepsilon }\setminus B)\to \pi_1(\mathbb B_{r_0}\setminus B)$, induced by ball imbedding, is an isomorphism.
\end{thm}
The group $\pi_1^{loc} (B,o):=\pi_1 (\mathbb B_{ r_0}\setminus B)$ is called the {\it local fundamental group} of $B$.

We say that an equisingular deformation $(\mathbb B_{\varepsilon}\times  D_{\delta},\mathcal B, \text{pr}_2)$ is {\it strong} if $\varepsilon< r(B_{\tau})$ for all $\tau\in D_{\delta}$. Let $l:[0,1]=\{ 0\leq t\leq 1\}\to D_{\delta}$ be a smooth path in $D_{\delta}$ and $(\mathbb B_{\varepsilon}\times  D_{\delta},\mathcal B, \text{pr}_2)$ a strong equisingular deformation.
In this case it is easy to show that
$\text{pr}_2: (\mathbb B_{\varepsilon}\times  D_{\delta}\setminus \mathcal B)\times_{D_{\delta}} [0,1] \to [0,1]$
is a $C^{\infty}$-trivial fibration with a fibre $\mathbb B_{\varepsilon}\setminus B_{l(0)}$ and, in particular, the fundamental groups $\pi_1((\mathbb B_{\varepsilon}\times  D_{\delta}\setminus \mathcal B)\times_{D_{\delta}} [0,1],(p_t,l(t))$ and
$\pi_1(\mathbb B_{\varepsilon}\setminus  B_{l(t)},p_{t})$ are naturally isomorphic for each $t\in [0,1]$ and $p_{t}\in \mathbb B_{\varepsilon}\setminus B_{l(t)}$.
Below, for each smooth path $l$ in $D_{\delta}$, we fix one of $C^{\infty}$-trivializations  of the fibration $\text{pr}_2: (\mathbb B_{\varepsilon}\times D_{\delta}\setminus \mathcal B)\times_{D_{\delta}} [0,1] \to [0,1]$ and call it an {\it equipment} of  $l$.

\subsection{$D$-automorphisms.} \label{def-act} Denote by $\mathfrak{B}_T$ the infinite dimensional subvariety of the variety of convergent power series in $\mathbb C[[z,w]]$, consisting of all power series $h(u,v)$ such that the germs of curves given by equations $h(u,v)=0$ in $(V,o)$, have the same type of singularity $T$ (i.e., all of these germs are equisingular deformation equivalent).
Consider a variety
$$\widetilde{\mathfrak{B}}_T =\{ \widetilde h=(h,\varepsilon,p)\in \mathfrak{B}_T\times\mathbb R\times\mathbb C^2 \mid h\in \mathfrak{B}_T, \varepsilon< r(B_h), p\in \mathbb B_{\varepsilon}\setminus B_{h}\}$$
and define in this variety smooth paths $l:[0,1]=\{ 0\leq t\leq 1\}\to \widetilde{\mathfrak{B}}_T$ of three types (below, {\it elementary admissible paths}); the smoothness of a path $l=(h_t,\varepsilon_t,p_t)$ means that the coefficients of the series $h_t$ and the points  $\varepsilon_t$ and $p_t$ smoothly depend on $t$). The first type consists of the paths  $l(t)= (h_t,\varepsilon_t,p_t)$ in which $h_t=h_0$ and $p_t=p_0$ for all $t\in [0,1]$. The second type consists of the paths  $l(t)= (h_t,\varepsilon_t,p_t)$ in which  $h_t=h_0$ and $\varepsilon_t=\varepsilon_0$ for all $t\in [0,1]$.
We call the paths of the first and second types the {\it $0$-paths}. The third type consists of the paths  (we call them the {\it $d$-paths}) $l(t)= (h_t,\varepsilon_t,p_t)$ for which
\begin{itemize} \item[$(i)$] $\varepsilon_t=\varepsilon_0$  for all $t\in [0,1]$,
\item[$(ii)$] there exist a smooth path $\widetilde l:[0,1]\to D_{\delta}$ and a strong equisingular deformation
$(\mathbb B_{\varepsilon}\times D_{\delta}, \mathcal B,\text{pr}_2)$ of the curve germ $B_{h_0}$ such that
$B_{h_t}=\mathcal B\cap \text{pr}_2^{-1}(\widetilde l(t))$  for all $t\in [0,1]$ and $\{ (p_t,t)\mid t\in [0,1]\}$ is a constant section of the equipment of $\widetilde l$.
\end{itemize}
Obviously, each elementary admissible path $l$ defines an isomorphism from the group $\pi_1(\mathbb B_{\varepsilon_{l(0)}}\setminus B_{h_{l(0)}},p_{l(0)})$ to $\pi_1(\mathbb B_{\varepsilon_{l(1)}}\setminus B_{h_{l(1)}},p_{l(1)})$.

{\it An admissible path} $l$ in $\widetilde{\mathfrak{B}}_T$ is a finite sequence of elementary admissible paths  $(l_1,\dots, l_n)$, in which the end of the path $l_i$ coincides with the beginning of the path $l_{i+1}$ for $1\leq i\leq n-1$.

We fix a point $\widetilde h_{0}=(h_0(u,v),\varepsilon_0,p_0)\in \widetilde{\mathfrak{B}}_T$ and call it  (and, resp., the curve germ $(B_{h_0},o)$ given by $h_0(u,v)=0$) a {\it base representative of singularity type} $T$.
Denote by $\Omega_T(\widetilde h_{0})$  the space of all admissible loops in $\widetilde{\mathfrak{B}}_T$ beginning at the point $\widetilde h_{0}$. Obviously, $\Omega_T(\widetilde h_{0})$ is a group and  a natural homomorphism
$$\text{{\bf Def}} :\Omega_T(\widetilde h_{0})\to \text{Aut}( \pi_1(\mathbb B_{\varepsilon_0}\setminus B_{h_0},p_0))=\text{Aut}(\pi_1^{loc}(B_{h_0},o))$$
is correctly defined.
The image  $\mathfrak{D}_T:=\text{{\bf Def}} (\Omega_T(\widetilde h_{0}))$ is called the {\it $D$-automorphism group.} Note that $\mathfrak{D}_T$ contains the group of inner automorphisms of the group $\pi_1^{loc}(B_{h_0},o)$.

\subsection{On the fundamental groups of the complements of curve germs.} \label{f-loc}
Zariski -- van Kampen Theorem (see below) gives an approach to finding representa\-tions of the local fundamental groups of curve germs. It is based on the calculation of the braid monodromy of the singularity and consists in the following (for a more detailed description, see, for example, \cite{K-U}).
Let a germ $(B,o)\subset (V,o)$ of a reduced holomorphic curve do not contain the germ $u=0$ and $m$
the multiplicity of the singularity of $B$ at the point $o$. Assume also that it is given by an equation
\begin{equation} \label{W}
v^m+\sum_{i=1}^{m}q_{i}(u)v^{m-i}=0,
\end{equation}
where $q_i(u)$ are convergent power series in $V=\mathbb B_{\varepsilon}$, $q_i(0)=0$, and
the polynomial $v^m+\sum q_{i}(u)v^{m-i}\in \mathbb C[[u]] [v]$ has not
multiple factors. Therefore, one can choose a small bidisk
$\overline D=\overline D_1\times \overline D_2\subset\mathbb B_{\varepsilon}$,
$$\overline D_1=\overline D_{1}(\varepsilon_1)=\{ u\in \mathbb C \, \mid \,\, \mid u\mid \leq \varepsilon_1 \, \}, \quad
\overline D_2=\overline D_{2}(\varepsilon_2)= \{ v\in \mathbb C \, \mid \,\, \mid v\mid \leq \varepsilon_2 \, \},$$
such that:
\begin{itemize}
\item[$(P1)$] the projection on the $u$-factor $\text{pr}=\text{pr}_1: B\cap \overline D\to \overline D_{1}$
is a proper finite map of degree $m$,
\item[$(P2)$] $|v|< \varepsilon_2$ for each point $(u,v)\in \text{pr}^{-1}(\overline D_1)$ and $o=(0,0)$ is the unique critical point of $\pr_{\mid B\cap\overline D}$.
\end{itemize}

Let us pick a point $(p_1,p_2)\in\partial \overline D_1\times \partial \overline D_2$, $p_1=\varepsilon_1$, $p_2=e^{3\pi i/2}\varepsilon_2$, and put
$$K_B:=K_B(\varepsilon_1)=\text{pr}_2(\text{pr}^{-1}(p_1))=\{v_1,\dots,v_m\}\subset \overline D_{2}.$$
The fundamental group $\pi_1 (\overline D_2\setminus K_B,p_2)\simeq \mathbb F_m$ is the free group of rank $m$ and it is generated by $m$ bypasses $\gamma_1,\dots,\gamma_m$ around the points $v_1\dots, v_m$. An ordered set
$ \{\gamma_1,\dots, \gamma_m\}$ is called a {\it  good geometric base} of $\pi_1(\overline D_2\setminus K_B, p_2)$ if the product
$\gamma_1\cdot .\, .\, .\, \cdot \gamma_m$ is equal to the element in $\pi_1 (\overline D_2\setminus K_B,p_2)$ represented by the circuit along the circle $\partial D_2$ in positive direction.

Note that $\text{pr}^{-1}( (\varepsilon_1t,\varepsilon_2))$, $0< t\leq 1$, is the disjoint union of $m$ paths. Therefore, hereafter, we can identify groups $\pi_1 (\overline D_2\setminus K_B(\varepsilon_1t),p_2)$ for different $t$ and if we need to decrease $\varepsilon_2$, we will identify the groups $\pi_1 (\overline D_2\setminus K_B,p_2)$ using the path $v(t)=e^{3\pi i/2}\varepsilon_2t$.

An example of a good geometric base of $\pi_1(\overline D_2\setminus K_B,p_2)$, where $B$ is given by equation $v((v-ku^k)^2-u^{2k+n})=0$, 
is shown in Fig. 1 in which $\varepsilon<\varepsilon_2$ and
$v_1= k\varepsilon^k+\sqrt{\varepsilon^{2k+n}}$,  $v_2=k\varepsilon^k-\sqrt{\varepsilon^{2k+n}}$, $v_3=0$ if $k\geq 1$, and
$v_1= \sqrt{\varepsilon^{n}}$, $v_2=0$, $v_3=-\sqrt{\varepsilon^{n}}$ if $k=0$.

\vspace{1cm}

\begin{picture}(400,100)

\put(185,100){\vector(0,-1){1}}
\put(200,100){\circle*{2}}
\put(203,100){$v_3$}
\put(200,100){\circle{30}}
\put(200,0){\circle*{2}}
\put(200,0){\line(0,1){85}}
\put(188,50){$\gamma_3$}
\put(188,0){$p_2$}
\put(265,100){\vector(0,-1){1}}
\put(280,100){\circle*{2}}
\put(283,100){$v_2$}
\put(280,100){\circle{30}}
\put(200,0){\line(3,4){68}}
\put(225,50){$\gamma_2$}
\put(305,100){\vector(0,-1){1}}
\put(320,100){\circle*{2}}
\put(323,100){$v_1$}
\put(320,100){\circle{30}}
\put(200,0){\line(4,3){113}}
\put(280,50){$\gamma_1$}
\put(190,-30){$\text{Fig.}\, 1$}
\end{picture}\vspace{2cm}

To recall the definition of the braid monodromy $b_{(B,o)}$ of $B$, consider  the braid group $\text{Br}_m$.
The group $\text{Br}_m$ has the following presentation. It generated by elements $a_1,\dots ,a_{m-1}$ being subject to relations
\begin{equation} \label{eqbr}
\begin{array}{clc} a_ia_{i+1}a_i & = & a_{i+1}a_i
a_{i+1} \qquad \qquad
1\leq i\leq m-1 ,  \\
a_ia_{k} & = & a_{k}a_i  \qquad \qquad \qquad \, \, \mid i-k\mid \,
\geq 2
\end{array}
\end{equation}
(such generators of the braid group $\text{Br}_m$ are called  {\it Artin's or standard}). The element $\Delta_m = (a_1a_2\dots a_{m-1})^m$ is called the {\it full twist}  and it belongs to the center of $\text{Br}_m$.

The loop $\partial \overline D_1$ oriented counter-clockwise and starting at $p_1$ lifts, via
$\text{pr}^{-1}(\partial \overline D_1)\cap B$, to
a motion $\pr_2(\{v_1(t),\dots,v_m(t)\})$ of $m$ distinct points in $\overline D_2$ starting and ending at $K_B$.
This motion defines a braid $b_{(B,o)}\in \text{Br}_m$ which is called the {\it braid monodromy of
$(B,o)$ with respect to} $\text{pr}$. 

\begin{lem} {\rm (\cite{K-U}, Lemma 4.1)}\label{cl1} Let $(B,o)$ be the singularity of type $A_{n}$, then
$$b_{(B,o)}= a_1^{n+1}.$$
\end{lem}

\begin{lem}\label{cl3} Let
$(B,o)$ be the singularity of type $T_{3,4k+n, 4k}$, then
$b_{(B,o)}= \Delta_3^ka_1^{n+1}.$
\end{lem}
\proof The braid $b_{(B,o)}$ consists of three threads
$$v_1(t)= \varepsilon^{2k}e^{4kt\pi i}+\sqrt{\varepsilon^{4k+n+1}}e^{(2k+\frac{n+1}{2})\pi i},
v_2(t)= \varepsilon^{2k}e^{4kt\pi i}-\sqrt{\varepsilon^{4k+n+1}}e^{(2k+\frac{n+1}{2})\pi i}, v_3(t)=0$$
in the space $\overline D_2\times \{ 0\leq t\leq 1\}$ starting at points of $K_B\times \{ 0\}=\{v_1,v_2,v_3\}\times \{ 0\}$; threads $v_1(t)$ and $v_2(t)$ make parallel $2k$ full twists around the thread $v_3(t)$ and make twists around each other $2k+\frac{n+1}{2}$ times. It is easy to see that the braid consisting of one parallel turn of $v_1(t)$ and $v_2(t)$ around $v_3(t)$ equals to $a_2a_1^2a_2$ and one half-twist of $v_1(t)$ and $v_2(t)$ is the braid $a_1$. It is easy to check that $a_2a_1^2a_2$ and $a_1$ commute in $\text{Br}_3$ and $(a_2a_1^2a_2)a_1^2=\Delta_3$. \qed

\begin{lem} {\rm (\cite{K-U}, Lemma 4.1)}\label{cl2} Let $(B,o)$ be the singularity of type $T_{3,n-1,n}$, then
$$b_{(B,o)}= (a_1a_2a_1)^{n}.$$
\end{lem}
\proof The braid $b_{(B,o)}$ consists of three threads
$$v_1(t)= \sqrt{\varepsilon^{n}}e^{\pi n i},\,\,\, v_2(t)=0, \,\,\,
v_3(t)= -\sqrt{\varepsilon^{n}}e^{\pi ni}$$
in the space $\overline D_2\times \{ 0\leq t\leq 1\}$ starting at points of $K_B\times \{ 0\}=\{v_1,v_2,v_3\}\times \{ 0\}$. Threads $v_1(t)$ and $v_3(t)$ make $n$ half-twists around the thread $v_2(t)$. It is easy to see that the braid consisting of one half-twist of $v_1(t)$ and $v_3(t)$ around $v_2(t)$ equals to $a_1a_2a_1=a_2a_1a_2$. \qed

\begin{lem}\label{Dn} Let $(B,o)$ be the singularity of type $D_{n}$, $n\geq 4$, then
$b_{(B,o)}= \Delta_3a_2^{n-4}.$
\end{lem}
\proof As a representative $(B,o)$ of the singularity type $D_{n}$, take the germ given by equation $(u-v)(v^2-u^{n-2})=0$. Then, as with the proof of the Lemma \ref{cl3}, it is easy to see that the braid $b_{(B,o)}$ consists of three strands, the first of which makes a full twist around the other two strands, and these two strands are twisted together $n/2-1$ times. \qed

\begin{lem} {\rm (\cite{K-U}, Lemma 4.1)}\label{clxx} Let $(B,o)$ be the singularity of type $T_{3}(n,\beta)$, then
$$b_{(B,o)}= \Delta_3^{n}(a_1a_2)^{\beta}.$$
\end{lem}

Fix a good geometric base $\gamma_1,\dots,\gamma_m$ of $\pi_1(\overline D_2\setminus K_B,p_2)$. The braid group $\text{Br}_m$ acts on the free group $\pi_1 (\overline D_2\setminus K_B,p_2)\simeq \mathbb F_m$ as follows
\begin{equation} \label{actbr}
\begin{array}{lcl} (\gamma_j)a_i &= & \gamma_j \qquad \qquad
\text{if}\, \, j\neq
i,i+1, \\
(\gamma_i)a_i &= & \gamma_i\gamma_{i+1}\gamma_i^{-1}, \\
(\gamma_{i+1})a_i &= & \gamma_i,
\end{array}
\end{equation}

In particular,
$$ (\gamma_i)\Delta_m=(\gamma_1\dots\gamma_m)\gamma_i(\gamma_1\dots\gamma_m)^{-1}\,\,\,\text{for}\,\,\,i=1,\dots ,m.$$

The imbedding $i:\{ p_1\}\times \overline D_2\hookrightarrow \overline D_1\times \overline D_2\hookrightarrow \mathbb B_{\varepsilon}$ defines a homomorphism
$$i_*:\pi_1(\overline D_2\setminus K_B,p_2)\to \pi_1(\mathbb B_{\varepsilon}\setminus B,p)\simeq \pi_1^{loc}(B,o).$$

\noindent {\bf Zariski -- van Kampen Theorem.} {\it The homomorphism $i_*:\pi_1(\overline D_2\setminus K_B,p_2)\to \pi_1^{loc}(B,o)$ is an epimorphism. The group $\pi_1^{loc}(B,o)$ has the following presentation}:
\begin{equation} \label{Zar-K} \pi_1^{loc}(B,o)=\langle \gamma_1,\dots ,\gamma_m \,\, \mid\,\, \gamma_i=(\gamma_i)b_{(B,o)},\,\, i=1,\dots, m\rangle , \end{equation}
where $\gamma_1,\dots,\gamma_m$ is a good geometric base of $\pi_1^{loc}(B,o)$. \\

The following two Lemmas are well known (see, for example, \cite{K-U}) and easily follows from Lemmas \ref{cl1} and \ref{clxx}.
\begin{lem} \label{cl1.2}  Let $(B,o)$ has a singularity of type $A_{n}$, where $n=2k-\delta$ and $\delta= 0$ or $1$. Then
\begin{equation}\label{A(n)}\pi_1^{loc}(B,o)=\langle \gamma_1, \gamma_2\, \, \mid \,\, (\gamma_1\gamma_2)^k\gamma_1^{1-\delta}=(\gamma_2\gamma_1)^k\gamma_2^{1-\delta}\rangle .\end{equation}
\end{lem}

\begin{lem} \label{clxx.1}  Let $(B_{n,\beta},o)$ has a singularity of type $T_3(n,\beta)$. If $\beta=1$ then
\begin{equation}\label{T(n,beta)} \begin{array}{ll} \pi_1^{loc}(B_{n,1},o)= \langle \gamma_1, \gamma_2,\gamma_3\, \, \mid & \gamma_1=(\gamma_1\gamma_2\gamma_3)^{n+1}\gamma_3(\gamma_1\gamma_2\gamma_3)^{-(n+1)}, \\
& \gamma_2=(\gamma_1\gamma_2\gamma_3)^{n}\gamma_1(\gamma_1\gamma_2\gamma_3)^{-n},
\\ & \gamma_3=(\gamma_1\gamma_2\gamma_3)^{n}\gamma_2(\gamma_1\gamma_2\gamma_3)^{-n}\rangle
\end{array} \end{equation}
and if $\beta=2$ then
\begin{equation}\label{T(n,beta2)} \begin{array}{ll} \pi_1^{loc}(B_{n,2},o)= \langle \gamma_1, \gamma_2,\gamma_3\, \, \mid & \gamma_1=(\gamma_1\gamma_2\gamma_3)^{n+1}\gamma_2(\gamma_1\gamma_2\gamma_3)^{-(n+1)}, \\
& \gamma_2=(\gamma_1\gamma_2\gamma_3)^{n+1}\gamma_3(\gamma_1\gamma_2\gamma_3)^{-(n+1)},
\\ & \gamma_3=(\gamma_1\gamma_2\gamma_3)^{n}\gamma_1(\gamma_1\gamma_2\gamma_3)^{-n}\rangle .
\end{array} \end{equation}
\end{lem}

\subsection{$WD$-subgroups of $D$-automorphism groups.} In this subsection, we use the notations and agreements of the previous paragraph. Denote by $\mathcal W_T$ a subspace of the space $\mathfrak B_T$ consisting of the power serious of the form (\ref{W}) and, respectively, $\widetilde{\mathcal W}_T=\{ \widetilde h=(h,\varepsilon, (\varepsilon_1, e^{3\pi i/2}\varepsilon_2))\in \widetilde{\mathfrak B}_T \mid h\in \mathcal W_T, \varepsilon_1,\varepsilon_2\in \mathbb R_+\}$, where $\varepsilon_1$ and $\varepsilon_2$ are such that   $h(u,v)$ has properties $(P1)$ and $(P2)$ in the bidisk $\overline D=\overline D_1(\varepsilon_1)\times \overline D_2(\varepsilon_2)\subset \mathbb B_{\varepsilon}$.

A $d$-path $l=(h_t,\varepsilon, (\varepsilon_1, e^{3\pi i/2}\varepsilon_2))$ in $\widetilde{\mathcal W}_T$ is called a $wd$-{\it path} if the restriction of the equipment of $l$ to
$\{ \partial_2 \overline D\times \{ l(t)\} \mid t\in  [0,1]\}$ is a trivial fibration
$\text{pr}_2:\partial_2 \overline D\times [0,1]\to [0,1]$,
where  $\partial_2 \overline D=\overline D_1(\varepsilon_1)\times \partial \overline D_2(\varepsilon_2)$.

A loop
$l=(l_1,\dots,l_k,l_{k+1},\dots, l_{k+n},l_k^{-1},\dots,l_1^{-1})$
in $\widetilde{\mathcal W}_T$, where $l_1,\dots,l_k$ are $0$-paths and $l_{k+1},\dots,l_{k+n}$ -- $wd$-paths, is called a
$w$-{\it loop}.

Fix a point $\widetilde h_{0}=(h_0(u,v),\varepsilon, (\varepsilon_1, e^{3\pi i/2}\varepsilon_2))$ in
$\widetilde{\mathcal W}_T\subset \widetilde{\mathfrak{B}}_T$ as the base representative of singularity type $T$ and denote by  $\Omega_{T,W}(\widetilde h_{0})$ the subgroup of
$\Omega_T(\widetilde h_{0})$ generated by the  $w$-loops in $\widetilde{\mathfrak{B}}_T$ beginning at the point $\widetilde h_{0}$.

Fix a good geometric base ofI $\pi_1(\overline D_2(\varepsilon_2)\setminus K_{B_{h_0}},e^{3\pi i/2}\varepsilon_2)$. Due to the identification of groups $\pi_1(\overline D_2(\lambda_2e^{3\pi i/2}\varepsilon_2)\setminus K_{B_r}(\lambda_1\varepsilon_1),\lambda e^{3\pi i/2}\varepsilon_2))$ for all $\lambda_1$, $\lambda_2\in (0,1)$ defined in Subsection \ref{f-loc}, the movements along
$w$-loops define a homomorphism
$\text{{\bf def}} :\Omega_{T,W} (\widetilde h_0)\to \text{Br}_m\subset\text{Aut}(\pi_1(\overline D_2(\varepsilon_2)\setminus K_{B_{h_0}},e^{3\pi i/2}\varepsilon_2))$
such that
$$ (i_*(\gamma))\text{{\bf Def}}(l)=i_*((\gamma)\text{{\bf def}}(l))$$ for all $l\in \Omega_{T,W} (\widetilde h_0)$  and $\gamma\in \pi_1(\overline D_2\setminus K_{B_{h_0}},p_2)$.
The images  $\mathfrak{d}_{T,W}:=\text{{\bf def}} (\Omega_{T,W} (\widetilde h_0))$ and $\mathfrak{D}_{T,W}:=\text{{\bf Def}} (\Omega_{T,W} (\widetilde h_0))\subset \mathfrak{D}_{T}$ are called the {\it  $WD$-automorphism groups.}

\begin{lem}\label{cl31} If $v((v-u^{k})^2-u^{2k+n+1})=0$ is the equation of $(B_{h_0},o)$, $k\geq 1$, $n\geq 0$, then the braids  $a_1$ and $a_2a_1^2a_2$ contain in $\mathfrak{d}_{T_{3,4k+n,4k},W}\subset Br_3$.
\end{lem}
\proof It is easy to see that the loop $l=\{ h_t \mid 0\leq t\leq 1\} $, given by
$$v[(v-u^{k})^2-e^{2\pi ti}u^{2k+n+1}]=0$$ in $\mathfrak{B}_{T_{3,4k+n,4k}}$ can be lifted in $\Omega_{T_{3,2k+n,2k},W}(\widetilde h_0)$ and $\text{\bf def}(\widetilde l)$ is the standard generator $a_1\in Br_3$. Similarly, the loop $l=\{ h_t\mid 0\leq t\leq 1\}$ given by
$v[(v-e^{2\pi ti}u^{k})^2-u^{2k+n+1}]=0$ in $\mathfrak{B}_{T_{3,2k+n,2k}}$ gives the element $a_2a_1^2a_2\in \mathfrak{d}_{T_{3,2k+n,2k}W}$. \qed

\begin{lem} \label{cl2.3} If $v(v^2-u^{2n})=0$ is the equation of $(B_{h_0},o)$, then $\mathfrak{d}_{T_{3,2n-1,2n},W}=Br_3$.
\end{lem}
\proof It is easy to see that for uplifts $\widetilde l_i$ in $\Omega_{T_{3,2n-1,2n},W}(\widetilde h_0)$ of the loops $l_i$, given by
$$(v+u^n)(v-\frac{1}{2}(1-e^{\pi ti})u^n)(v-\frac{1}{2}(1+e^{\pi ti})u^n)=0,$$
$$(v+\frac{1}{2}(1+e^{\pi ti})u^n)(v+\frac{1}{2}(1-e^{\pi ti})u^n)(v-u^n)=0$$
in $\mathfrak{B}_{T_{3,2n-1,2n}}$, their images $\text{\bf def}(\widetilde l_i)$ are the standard generators $a_1$ and $a_2$ of the braid group $Br_3$ acting on
$\pi_1(\overline D_2\setminus K_{B_{h_0}})$.   \qed

\subsection{ Deformation equivalence of germs of 
covers.}\label{defm}
Consider the germ of a finite cover $F: (U,o')\to (V,o)$, $\deg F=d$. Choose local coordinates $z,w$ in $(U,o')$ and $u,v$ in $(V,o)$. The germ $F$ is given by two functions
$$ \begin{array}{ll}
 u= & f_1(z,w), \\
 v= & f_2(z,w), \end{array}
$$
where $f_i(z,w)\in H^0(U,\mathcal O_U)$. The ramification divisor $R$ in $U$ is defined by equation
$$J(f):= \det \left(\begin{array}{cc} \frac{\partial u}{\partial z} & \frac{\partial u}{\partial w} \\ \frac{\partial v}{\partial z} & \frac{\partial v}{\partial w}\end{array}\right)= 0 $$
and let $B=F(R_{red})\subset (V,o)$ be the branch curve of the germ of finite cover $F$.
\begin{rem} \label{remark0} The divisor $R\subset (U,o')$ and the curve germ $B\subset (V,o)$ depend only on $F$ and do not depend on the choice of coordinates in $(U,o')$ and $(V,o)$. \end{rem}

\begin{df} \label{def-d} Let $\mathcal F: (U,o')\times D_{\delta} \to (V,o)\times D_{\delta}$ be a finite holomorphic mapping branched along a surface $\mathcal B\subset (V,o)\times D_{\delta}$ and
such that $\text{pr}_2\circ \mathcal F=\text{Pr}_2$, where  $\text{Pr}_2: (U,o')\times D_{\delta} \to D_{\delta}$ and $\text{pr}_2: (V,o)\times D_{\delta} \to D_{\delta}$ are the projections to the second factor. The cover $\mathcal F$ is called a strong deformation of the germ of finite cover $F_0=\mathcal F_{\mid (U,o')\times \{ 0\}}: (U,o')\times \{ 0\}\to (V,o)\times \{ 0\}$ if
$((V,o)\times D_{\delta}, \mathcal B, 
\text{pr}_2)$ is a strong equisingular deformation of the curve germ $B_0=\mathcal B\cap \text{pr}_2^{-1}(0)$ and the germs of finite covers  $F_{\tau}=\mathcal F_{\mid (U,o')\times \{ \tau \}}: (U,o')\times \{ t\}\to (V,o)\times \{ \tau \}$ are called strong  deformation equivalent to the germ $F_0$.
\end{df}

Let us continue the  strong deformation equivalence of germs of finite covers to an equivalence relation.

The germ of finite cover $F:(U,o')\to (V,o)$ defines a homomorphism
$$F_*:\pi_1^{loc}(B,o)=\pi_1(V\setminus B,p)\to \mathbb S_d,$$
where $\mathbb S_d$ is the symmetric group acting on the fibre $F^{-1}(p)$. Note that the homomorphism $F_*$ is defined uniquely only if we fix a numbering of the points of $F^{-1}(p)$ and in general case it is defined uniquely up to inner automorphism of $\mathbb S_d$. The group $G_F=\text{im} F_*$ is called the {\it (local) monodromy group} of the germ $F$. The group $G_F$ is a transitive subgroup of $\mathbb S_d$, since $(U,o')$ is an irreducible germ of a smooth surface and hence $U\setminus F^{-1}(B)$ is connected.
By Grauert - Remmert - Riemann - Stein Theorem (\cite{St}), the epimorphism $F_*:\pi_1(V\setminus B)\to G_F\subset \mathbb S_d$ uniquely determines the germ of a finite cover $F$.

Let $(V\times D_{\delta}, \mathcal B, \text{pr}_2)$ \hspace{0.1cm} be a strong equisingular deformation of a germ $(B,o)=\mathcal B\cap\text{pr}_2^{-1}(0)$. 
By  Grauert - Remmert - Riemann - Stein Theorem, the homomorphism $F_*: \pi_1^{loc}(B,o)\simeq \pi_1((V\times D_{\delta})\setminus\mathcal B)\to \mathbb S_d$ defines a finite $d$-sheeted cover $\mathcal F: \mathcal U\to V\times D_{\delta}$, where in general case $\mathcal U$ is a normal complex-analytic variety. But, if $F_*$ is defined by a germ of finite cover $F:(U,o')\to (V,o)$ branched in $(B,o)$ then we have
\begin{prop}\label{eqi-cov} The cover $\mathcal F: \mathcal U\to V\times D_{\delta}$ is a strong deformation of the germ of finite cover $F=\mathcal F_{\mid \mathcal F^{-1}((V,o)\times \{ 0\})}: (U,o')\times \{ 0\}\to (V,o)\times \{ 0\}$.
\end{prop}
\proof First  of all, note that $\mathcal F^{-1}(V\times \{0\})=(U,o')$ and $\mathcal F_{\mid U}=F$.

In notations of Definition \ref{def1}, the homomorphism $F_*$ defines a finite cover
$\widetilde{\mathcal F}: \widetilde{\mathcal U}\to\mathcal V_n$ branched in $\mathcal B_n$. Since $\mathcal B_n$ is a divisor with normal crossings and all its singular points are sections of $\text{pr}_2\circ\widetilde{\sigma}_1\circ\dots\circ\widetilde{\sigma}_{n}$, then the local fundamental groups of the complement to $\mathcal B_n$ at the points of $\mathcal B_n$ are abelian and hence  $Sing\, \widetilde{\mathcal U}$ is, first,  a disjoint union $\bigsqcup \widetilde S_j$ of sections $\widetilde S_j$ of $\text{pr}_2\circ\widetilde{\sigma}_1\circ\dots\circ\widetilde{\sigma}_{n}\circ \widetilde F$ lying over $Sing\, \mathcal B_n$ and second, at the points of $\widetilde S_{j_0}\subset Sing\, \widetilde{\mathcal U}$, the variety $\widetilde{\mathcal U}$ locally biholomorphic to $W_{j_0}\times D_{\delta_1}$, where $W_{j_0}$ is a germ of two-dimensional cyclic quotient singularity depending on the local monodromy at  the points of $S_{j_0}\subset Sing\, \mathcal B_n$ (see details in \cite{B}, III.6). The minimal resolutions of singularities $\rho_{j} : \overline W_j\to W_j$ defines a resolution of singularities $\rho:\overline{\mathcal U}\to\widetilde{\mathcal U}$. By Stein Factorisation Theorem, there is a holomorphic mapping $\Sigma: \overline{\mathcal U}\to \mathcal U$ contracting the divisor $(\mathcal F\circ\Sigma)^{-1}(\mathcal E)$ to the section $\mathcal F^{-1}(\{ o\}\times D_{\varepsilon})$ of projection $\mathcal F\circ \text{pr}_2$, where $\mathcal E$ is the exceptional divisor of the sequence of the monoidal transformations  $\widetilde{\sigma}_1\circ\dots\circ\widetilde{\sigma}_n$.

Consider the restriction $\sigma:= \Sigma_{\mid\Sigma^{-1}((U,o'))}: \overline U=\Sigma^{-1}((U,o'))\to (U,o')$. The exceptional divisor $E$ of $\sigma$ is $E=\mathcal E\cap \overline U$. By Zariski Theorem $\sigma$ is a composition of $\sigma$-processes, since $\overline U$ and $U$ are nonsingular and $\sigma$ is a bimeromorphic holomorphic mapping. It
follows from \cite{Har} (see Chapter 2, Example 6.2.2.) and Nokano Contractibility Criterion (\cite{No}) that $\Sigma$ is the composition of monoidal transformations (being in one-to-one correspondence with the composition of $\sigma$-processes $\sigma$) of smooth threefolds with centers in sections of the projection to $D_{\delta}$. \qed \\

The following proposition easily follows from the proof of Proposition \ref{eqi-cov}.

\begin{prop}\label{equiD} If $F_1:(U,o')\to (V,o)$ and $F_2:(U,o')\to (V,o)$ are deformation equivalent covers, then their ramification divisors $R_{1,red}$ and $R_{2,red}$ are equisingular deformation equivalent.
\end{prop}

Let $\gamma_1,\dots, \gamma_m$ be a good geometric base of the fundamental group $\pi_1^{loc}(B,o)$ and $F_*:\pi_1^{loc}(B,o)\to \mathbb S_d$ a homomorphism to the symmetric group $\mathbb S_d$ such that $G_F=\text{im}\, F_*$ is a transitive subgroup of $\mathbb S_d$. The collection $\{ C_1,\dots ,C_m\}$ of conjugacy classes in $\mathbb S_d$, $F_*(\gamma_i)\in C_i$ for $i=1,\dots,m$, is called the {\it dataset} of the homomorphism $F_*$. We say that the homomorphism $F_*$ is {\it sole} if it is  uniquely, up to inner automorphisms of $\mathbb S_d$, defined by its dataset.

We say that two homomorphisms $F_{i*}:\pi_1^{loc}(B,o)\to\mathbb S_d$, $i=1,2$, are {\it equivalent} if they differ from each other by an inner automorphism of $\mathbb S_d$ and they are {\it deformation equivalent} if they differ from each other on a $D$-automorphism of $\mathfrak{D}_T$ and an inner automorphism of $\mathbb S_d$.

From the above it follows
\begin{cor}\label{eqi-cov1} Two germs of finite covers $F_1: (U,o')\to (V,o)$ and $F_2: (U,o')\to (V,o)$ of degree $d$, branched along $(B_1,o)\subset (V,o)$ and $(B_2,o)\subset (V,o)$, are deformation equivalent if and only if the curve germs $(B_1,o)$ and $(B_2,o)$ are equisingular deformation equivalent and the monodromies $F_{1*}: \pi_1\to \mathbb S_d$ and $F_{2*}:\pi_1\to \mathbb S_d$ are deformation equivalent, where $\pi_1=\pi_1^{loc}(B_1,o)=\pi_1^{loc}(B_2,o)$.
\end{cor}

\begin{cor}\label{eqi-cov2}  Two  deformation equivalent germs of finite covers $F_1: (U,o')\to (V,o)$ and
$F_2: (U,o')\to (V,o)$ of degree $d$, branched along $(B_1,o)\subset (V,o)$ and $(B_2,o)\subset (V,o)$, are equivalent if the curve germ $(B_1,o)$ is rigid and the monodromy $F_{1*}: \pi_1^{loc}(B_1,o)\to \mathbb S_d$ is sole.
\end{cor}

\begin{lem}\label{cl7} Let $T_{3,8k+2n+4, 8k+4}$ be the singularity type of $(B,o)$
and a homomorphism $F_*:\pi_1^{loc}(B,o)\to \mathbb S_4$  have the following properties:
\begin{itemize}
\item[$(i)$] the image $F_*(\pi_1^{loc}(B,o))$ is a transitive subgroup of $\mathbb S_4$,
\item[$(ii)$] $F_*(\gamma_i)=\tau_i$ are transpositions, where $\gamma_1,\gamma_2,\gamma_3$ is a good geometric base.
\end{itemize}
Then $F_*$ is sole and $F_*(\pi_1^{loc}(B_r,o))=\mathbb S_4$.
\end{lem}
\proof By Lemma \ref{cl3} and Zariski -- van Kampen Theorem, the group $\pi_1^{loc}(B,o)$ has the following presentation
\begin{equation} \label{pres1}\begin{array}{ll}  \pi_1^{loc}(B,o)=\langle \gamma_1,\gamma_2,\gamma_3 \mid & \gamma_1=(\gamma_1\gamma_2\gamma_3)^{4k+2}(\gamma_1\gamma_2)^{n}\gamma_1\gamma_2\gamma_1^{-1}(\gamma_1\gamma_2)^{-n}
(\gamma_1\gamma_2\gamma_3)^{-(4k+2)},\\
& \gamma_2=(\gamma_1\gamma_2\gamma_3)^{4k+2}(\gamma_1\gamma_2)^{n}\gamma_1(\gamma_1\gamma_2)^{-n}
(\gamma_1\gamma_2\gamma_3)^{-(4k+2)}, \\
& \gamma_3=(\gamma_1\gamma_2\gamma_3)^{4k+2}\gamma_3(\gamma_1\gamma_2\gamma_3)^{-(4k+2)}  \rangle
\end{array} .\end{equation}
Since $F_*(\pi_1^{loc}(B,o))$ is a transitive subgroup of $\mathbb S_4$ and $F_*(\gamma_i)=\tau_i$ are transpositions, then $F_*(\pi_1^{loc}(B,o))=\mathbb S_4$ and $\tau_1\tau_2\tau_3$ is a cycle of length $4$. Therefore, with up to conjugation in $\mathbb S_4$, we can assume that $\tau_1\tau_2\tau_3=(1,2,3,4)$ and $\tau_3=(1,3)$, since $\gamma_3=(\gamma_1\gamma_2\gamma_3)^{4k+2}\gamma_3(\gamma_1\gamma_2\gamma_3)^{-(4k+2)}$. Consequently, $\tau_1\tau_2=[((1,2)(3,4)]$ and we can assume that $\tau_1=(1,2)$ and $\tau_2=(3,4)$ (if $\tau_1=(3,4)$ and $\tau_2=(1,2)$ then we conjugate by $[(1,3)(2,4)]$).\qed

\begin{lem}\label{cl8} Let
$(B,o)$ have the singularity type $T_{3,8m+2n, 8m}$ and let $F_*:\pi_1^{loc}(B,o)\to \mathbb S_4$ be a homomorphism such that
\begin{itemize}
\item[$(i)$] $F_*(\pi_1^{loc}(B,o))$ is a transitive subgroup of $\mathbb S_4$,
\item[$(ii)$] $F_*(\gamma_i)=\tau_i$ are transpositions, where $\gamma_1,\gamma_2,\gamma_3$ is a good geometric base.
\end{itemize}
Then $n=3k+1$, $k\in \mathbb Z_{\geq 0}$,  $F_*(\pi_1^{loc}(B,o))=\mathbb S_4$, and $F_*$ is defined uniquely with up to inner automorphisms of $\mathbb S_4$ and $D$-automorphisms of $\mathfrak{D}_{T_{3,8m+6k+2, 8m}}$. \end{lem}
\proof By Lemma \ref{cl3} and Zariski -- van Kampen Theorem, the group $\pi_1^{loc}(B,o)$ has the following presentation
\begin{equation} \label{pres2}\begin{array}{ll}  \pi_1^{loc}(B,o)=\langle \gamma_1,\gamma_2,\gamma_3 \mid & \gamma_1=(\gamma_1\gamma_2\gamma_3)^{4m}(\gamma_1\gamma_2)^{n}\gamma_1\gamma_2\gamma_1^{-1}(\gamma_1\gamma_2)^{-n}
(\gamma_1\gamma_2\gamma_3)^{-4m},\\
& \gamma_2=(\gamma_1\gamma_2\gamma_3)^{4m}(\gamma_1\gamma_2)^{n}\gamma_1(\gamma_1\gamma_2)^{-n}
(\gamma_1\gamma_2\gamma_3)^{-4m}, \\
& \gamma_3=(\gamma_1\gamma_2\gamma_3)^{4m}\gamma_3(\gamma_1\gamma_2\gamma_3)^{-4m} \,\, \rangle .
\end{array} \end{equation}
As in the proof of Lemma \ref{cl7}, since $F_*(\pi_1^{loc}(B,o))$ is a transitive subgroup of $\mathbb S_4$ and $F_*(\gamma_i)=\tau_i$ are transpositions, then $F_*(\pi_1^{loc}(B,o))=\mathbb S_4$ and $\tau_1\tau_2\tau_3$ is a cycle of length $4$ and $\tau_1\neq \tau_2$. It follows from (\ref{pres2}) that
\begin{equation} \label{rel1}
\tau_1=(\tau_1\tau_2)^{n}\tau_1\tau_2\tau_1^{-1}(\tau_1\tau_2)^{-n},\qquad \tau_2=(\tau_1\tau_2)^{n}\tau_1(\tau_1\tau_2)^{-n},
\end{equation}
and hence, $\tau_1$ and $\tau_2$ are not commute, since they are conjugated in the group $\langle \tau_1,\tau_2\rangle $.
Therefore, with up to conjugation in $\mathbb S_4$, we can assume that $\tau_1=(1,2)$, $\tau_2=(1,3)$, and $\tau_1\tau_2=(1,2,3)$. Now, it is easy to see that equalities (\ref{rel1}) hold iff $n=3k+1$.

There are three possibilities for $\tau_3$: either $\tau_3=(1,4)$ (denote this homomorphism by $F_1$), or $\tau_3=(2,4)$ (denote this homomorphism by $F_2$), or $\tau_3=(3,4)$ (denote this homomorphism by $F_3$).
Note that, by Lemma \ref{cl31}, $a_1\in \mathfrak{d}_{T_{3,8m+6k+2, 8m},W}$. Then the homomorphism  $F_{1*}$ and the homomorphism
$\widetilde F_{1*}$ sending $\gamma_1$ to $(1,3)$, $\gamma_2$ to $(2,3)$, and $\gamma_3$ to $(1,4)$ differ by the action of $D$-automorphism $i_*(a_1)$. It is easy to check that the homomorphism $\widetilde F_{1*}$ coincides with $F_{2*}$ after conjugation by $(1,2,3)$. Similarly, the homomorphism  $F_{2*}$ and the homomorphism $\widetilde F_{2*}$ sending $\gamma_1$ to $(1,3)$, $\gamma_2$ to $(2,3)$, and $\gamma_3$ to $(2,4)$ also differ by the action of $D$-automorphism $i_*(a_1)$. It is easy to check that the homomorphism $\widetilde F_{2*}$ coincides with $F_{3*}$ after conjugation by $(1,2,3)$.   \qed

\begin{lem}\label{cl181} Let
$(B,o)$ have the singularity type $T_{3,8m+2n-1, 8m}$, $m\geq 1$, $n\geq 0$, and $F_*:\pi_1^{loc}(B,o)\to \mathbb S_4$  a homomorphism such that
\begin{itemize}
\item[$(i)$]  $F_*(\pi_1^{loc}(B,o))$ is a transitive subgroup of $\mathbb S_4$,
\item[$(ii)$] $F_*(\gamma_i)=\tau_i$ are transpositions, where $\gamma_1,\gamma_2,\gamma_3$ is a good geometric base.
\end{itemize}
Then $F_*(\pi_1^{loc}(B,o))=\mathbb S_4$ and if either $n=3k+1$, or $n=3k+2$, $k\in \mathbb Z_{\geq 0}$, or $n=0$, then with up to deformation equivalence,
there exist the unique  such homomorphism  $F_*$, and if $n=3k>0$ then
there are at most two homomorphisms satisfying $(i)$ and $(ii)$.
\end{lem}
\proof By Lemma \ref{cl3} and Zariski -- van Kampen Theorem, the group $\pi_1^{loc}(B,o)$ has the following presentation
\begin{equation} \label{pres2}\begin{array}{ll}  \pi_1^{loc}(B,o)=\langle \gamma_1,\gamma_2,\gamma_3 \mid & \gamma_1=(\gamma_1\gamma_2\gamma_3)^{4m}(\gamma_1\gamma_2)^{n}\gamma_1(\gamma_1\gamma_2)^{-n}
(\gamma_1\gamma_2\gamma_3)^{-4m},\\
& \gamma_2=(\gamma_1\gamma_2\gamma_3)^{4m}(\gamma_1\gamma_2)^{n}\gamma_2(\gamma_1\gamma_2)^{-n}
(\gamma_1\gamma_2\gamma_3)^{-4m}, \\
& \gamma_3=(\gamma_1\gamma_2\gamma_3)^{4m}\gamma_3(\gamma_1\gamma_2\gamma_3)^{-4m} \,\, \rangle .
\end{array} \end{equation}
As in the proof of Lemma \ref{cl7}, since $F_*(\pi_1^{loc}(B,o))$ is a transitive subgroup of $\mathbb S_4$ and $F_*(\gamma_i)=\tau_i$ are transpositions, then $F_*(\pi_1^{loc}(B,o))=\mathbb S_4$. Therefore $\tau_1\neq \tau_2$ and $\tau_1\tau_2\tau_3$ is a cycle of length $4$.

It follows from (\ref{pres2}) that
\begin{equation} \label{rel1}
\tau_1=(\tau_1\tau_2)^{n}\tau_1(\tau_1\tau_2)^{-n},\qquad \tau_2=(\tau_1\tau_2)^{n}\tau_2(\tau_1\tau_2)^{-n}.
\end{equation}
Let $n=3k+1$ or $3k+2$. Then it follows from (\ref{rel1}) that $\tau_1$ and $\tau_2$ must commute with each other
and hence, with up to conjugation, we can assume that $\tau_1=(1,2)$ and $\tau_2=(3,4)$. Then $\tau_3$ belongs to the set
$\{ (1,3),(1,4), (2,3),2,4)\}$. But again, it is easy to see that there is an inner automorphism of $\mathbb S_4$ saving fixed
$\tau_1$ and $\tau_2$ and sending $\tau_3$ to $(1,3)$.

Since $\tau_1\tau_2\tau_3$ is a cycle of length $4$, we can assume that $\tau_1\tau_2\tau_3=(1,2,3,4)$ and Lemma \ref{cl181} in case $n=0$ follows from Theorem 2.1 in \cite{K1} and  Lemma \ref{cl2.3}.

Let $n=3k>0$. With up to conjugation in $\mathbb S_4$, we have only two possibilities: either $\tau_1=(1,2)$ and $\tau_2=(1,3)$, or $\tau_1=(1,2)$ and $\tau_2=(3,4)$. 
The case when $\tau_1=(1,2)$ and $\tau_2=(3,4)$ was considered above and the case when $\tau_1=(1,2)$ and $\tau_2=(1,3)$ was considered in the proof of Lemma \ref{cl8}. \qed

\begin{lem}\label{AnD4} Let $A_{4n-1}$ be the singularity type of a curve germ $(B=B_1\cup B_2,o)$ and a homomorphism
$F_*:\pi_1^{loc}(B,o)\to\mathbb S_4$ have the following properties:
\begin{itemize}
\item[$(i)$]  $G_F=\text{im}, F_*$ is a transitive group of $\mathbb S_4$;
\item[$(ii)$] $F_*(\gamma_1)=\tau_1\tau_2$ is the product of two commuting transpositions $\tau_1$ and $\tau_2$,  and $F_*(\gamma_2)=\tau_3$ is a transposition, where $\gamma_1$ is a bypass around $B_1$, $\gamma_2$ is a bypass around $B_2$, and $\gamma_1, \gamma_2$ is a good geometric base of $\pi_1^{loc}(B,o)$;
  \end{itemize}
Then $G_F=\mathbb D_4$ is a dihedral subgroup of $\mathbb S_4$ and $F_*$ is sole.
\end{lem}
\proof With up to conjugation in $\mathbb S_4$, we can assume that $\tau_1=(1,2)$ and $\tau _2=(3,4)$. Note that $\tau_3\neq \tau_i$ for $i=1,2$, since $G_F$ is a transitive subgroup of $\mathbb S_4$
and $\gamma_1$, $\gamma_2$ generate the group $\pi_1^{loc}(B,o)$. Therefore $\tau_3=(i_1,i_2)$, where $i_1\in \{ 1,2 \}$ and
$i_2\in \{ 3,4\}$. Again, applying a conjugation in $\mathbb S_4$, we can assume that $\tau_3=(1,3)$. \qed

\begin{lem} \label{Lem} Let the singularity type of a germ $(B=B_1\cup B_2,o)$ {\rm (}$v-2u=0$ is an equation of the germ $B_1$ and $v^2-u^{2k+1}=0$ is the equation of the germ $B_2${\rm )} be $D_{2k++3}$ and a homomorphism $F_*:\pi_1^{loc}(B,o)\to \mathbb S_4$  such that
\begin{itemize}
\item[$(i)$] $G_F=\text{im}\, F_*$ is a transitive subgroup of $\mathbb S_4$,
\item[$(ii)$] $F_*(\gamma_1)$ is a product of two different commuting transpositions and $F_*(\gamma_2)$, $F_*(\gamma_3)$ are transpositions,
\end{itemize}
 where $\gamma_1,\gamma_2,\gamma_3$ is a good geometric base of $\pi_1^{loc}(B,o)$, $\gamma_2$ and $\gamma_3$ are bypasses around $B_2$ and $\gamma_1$ is a bypass around $B_1$. Then $G_F$ is a dihedral subgroup $\mathbb D_4$ of $\mathbb S_4$ and $F_*$ is solo.
\end{lem}
\proof By Lemma \ref{Dn} and Zariski -- van Kampen Theorem, the group $\pi_1^{loc}(B,o)$ has the following presentation
\begin{equation} \label{pres5}\begin{array}{ll}  \pi_1^{loc}(B,o)=\langle \gamma_1,\gamma_2,\gamma_3 \mid & \gamma_1=(\gamma_1\gamma_2\gamma_3)\gamma_1(\gamma_1\gamma_2\gamma_3)^{-1}, \\ & \gamma_2=(\gamma_1\gamma_2\gamma_3)(\gamma_2\gamma_3)^{k}\gamma_2\gamma_3\gamma_2^{-1}(\gamma_2\gamma_3)^{-k}
(\gamma_1\gamma_2\gamma_3)^{-1},\\
& \gamma_3=(\gamma_1\gamma_2\gamma_3)(\gamma_2\gamma_3)^{k}\gamma_2(\gamma_1\gamma_2)^{-k}
(\gamma_1\gamma_2\gamma_3)^{-1}
 \rangle
\end{array} \end{equation}
Note that $\nu=F_*(\gamma_1)F_*(\gamma_2)F_*(\gamma_3)$ is an even permutation. Therefore $\nu$ is either a cycle of length $3$ or $\nu\in Kl_4$. But, it follows from (\ref{pres5}) that $\nu$ can not be a cycle of length $3$, since it should commute with $F_*(\gamma_1)$, and $\nu\neq 1$, since then $F_*(\gamma_2)F_*(\gamma_3)=F_*(\gamma_1)$ and $F_*(\gamma_1),F_*(\gamma_2),F_*(\gamma_3)$ can not generate a transitive subgroup of $\mathbb S_4$. Note also that  $F_*(\gamma_2)\neq F_*(\gamma_3)$, since, otherwise, it follows from (\ref{pres5}) that $F_*(\gamma_2)=F_*(\gamma_3)$ commutes with $F_*(\gamma_1)$ and hence, again,  $F_*(\gamma_1),F_*(\gamma_2),F_*(\gamma_3)$ can not generate a transitive subgroup of $\mathbb S_4$. Consequently,
$F_*(\gamma_2) F_*(\gamma_3)\in Kl_4$ and, with up to conjugation in $\mathbb S_4$, we can assume that $F_*(\gamma_1)=[(1,2)(3,4)]$,
$F_*(\gamma_2)=(1,3)$, and $F_*(\gamma_3)=(2,4)$. \qed

\begin{lem} \label{Lemx} Let the singularity type of a germ $(B,o)$ be $T_3(4n+\beta,\beta)$ and a homomorphism $F_*:\pi_1^{loc}(B,o)\to \mathbb S_4$ be such that
\begin{itemize}
\item[$(i)$] $G_F=\text{im}\, F_*$ is a transitive subgroup of $\mathbb S_4$,
\item[$(ii)$] $F_*(\gamma_i)$, $i=1,2,3$, are transpositions,
\end{itemize}
 where $\gamma_1,\gamma_2,\gamma_3$ is a good geometric base of $\pi_1^{loc}(B,o)$. Then $G_F=\mathbb S_4$ and $F_*$ is  solo.
\end{lem}
\proof
It follows from $(i)$ and $(ii)$ that $\nu=F_*(\gamma_1)F_*(\gamma_2)F_*(\gamma_3)$ is a cycle of length four. Without loss  of generality, we can assume that $\nu=(1,2,3,4)$.
Consider the case $\beta=1$ (the case $\beta=2$ is similar and therefore it will be omitted). By Lemma \ref{clxx.1}, we have
\begin{equation}\label{T(n,beta)1}  F_*(\gamma_1)=\nu^{2}F_*(\gamma_3)\nu^{-2},\quad
 F_*(\gamma_2)=\nu F_*(\gamma_1)\nu^{-1},\quad
 F_*(\gamma_3)=\nu F_*(\gamma_2)\nu^{-1}. \end{equation}
The element $\nu$ acts by conjugation on the set of transpositions in $\mathbb S_4$. There are two orbits of this action:
$$ (1,2)\mapsto (1,4)\mapsto (3,4)\mapsto (2,3)\mapsto (1,2),\,\,\, \text{and}\,\,\, (1,3)\mapsto(2,4)\mapsto (1,3). $$
Therefore, by $(i)$ and since we can conjugate on the element $\nu$, it follows from (\ref{T(n,beta)1}) that we can put
$F_*(\gamma_1)=(1,2)$, $F_*(\gamma_2)=(1,4)$, $F_*(\gamma_3)=(3,4)$, and after that to check that $\nu=[(1,2)]\cdot [(1,4)]\cdot [(3,4)]$. \qed

\section{Proof of Theorem \ref{main}}
Let $L_1=\{u_0=0\}$ and $L_2=\{v_0=0\}$ be the axes of some local complex-analytic coordinates $u_0,v_0$ in $(V,o)$. Then the local intersection number of the divisors $M_1=F^*(L_1)$ and $M_2=F^*(L_2)$ at the point $o'$ is equal to $(M_1,M_2)_{o'}=\deg_{o'} F=4$. Therefore we have the following possibilities:
\begin{itemize}
\item[$({\rm I})$] either $M_1$ or $M_2$ is a germ of a non-singular curve,
 \item[$({\rm II})$] $M_1$ and $M_2$ have a singularity of multiplicity 2 at the point $o'$.
\end{itemize}

Let $R\subset (U,o')$ be the ramification divisor of the germ of finite cover $F$ and $B=F(R_{red})\subset (V,o)$ the branch curve.

Denote $\mathfrak{m}\subset \mathbb C[[z_0,w_0]]$ the maximal ideal in the ring of power series $\mathbb C[[z_0,w_0]]$.

\subsection{Case $({\rm I})$.} \label{secti} Let $M_1$ be non-singular. Then we can choose local coordinates $z_0,w_0$ in $(U,o')$ such that $F^*(u_0)=z_0$ and $F^{*}(v_0)=v_0(z_0,w_0)=\sum_{i=0}^{\infty}a_i(z_0)w_0^i$, where
$a_i(z_0)=\sum_{j=0}^{\infty}a_{i,j}z_0^j\in \mathbb C[[z_0]]$. Performing the coordinates change $v_0 \leftrightarrow v_0-a_0(u_0)$, we can assume that $a_0(z_0)\equiv 0$. In addition, we have $a_{1,0}=a_{2,0}=a_{3,0}=0$ and can assume that $a_{4,0}=1$, since $(M_1,M_2)_{o'}=4$.

The divisor $R$ is given by equation
$$J(F):= \det \left(\begin{array}{cc} 1 & 0 \\ \frac{\partial v_0}{\partial z_0} & \frac{\partial v_0}{\partial w_0}\end{array}\right)=0, $$
i.e., $R$ is given by equation
\begin{equation}\label{eq1} \sum_{i=1}^{\infty}ia_i(z_0)w_0^{i-1}=0.\end{equation}
Let us write equation (\ref{eq1}) in the following form
\begin{equation} \label{eq2} z_0H_1(z_0,w_0)+4w_0^3+ H_2(z_0,w_0)=0,
\end{equation}
where $H_1(z_0,w_0)$ is a polynomial of degree $2$ and $H_2(z_0,w_0)\in \mathfrak{m}^4$. It follows from (\ref{eq2}) that $(M_1,R)_{o'}=3$. Therefore there are seven  possibilities:
\begin{itemize}
\item[$({\rm I}_1)$]  $R=3R_1$, where $R_1$ is a germ of a smooth curve and $(M_1,R_1)_{o'}=1$,
\item[$({\rm I}_2)$] $R=2R_1+R_2$, where $R_1$ and $R_2$ are germs of smooth curves $(M_1,R_1)_{o'}=(M_1,R_2)_{o'}=1$,
\item[$({\rm I}_3)$]  $R=R_1$ is reduced irreducible and $(M_1,R_1)_{o'}=3$,
\item[$({\rm I}_{4.1})$] $R=R_1+R_2$, where the germ $R_1$ is irreducible,
$(M_1,R_1)_{o'}=2$ and  $\deg f_{\mid R_1}=1$, $R_2$ is a germ of smooth curve and $(M_1,R_2)_{o'}=1$,
\item[$({\rm I}_{4.2})$] $R=R_1+R_2$, where the germ $R_1$ is irreducible,
$(M_1,R_1)_{o'}=2$, and $\deg f_{\mid R_1}=2$, $R_2$ is a germ of smooth curve and $(M_1,R_2)_{o'}=1$,
\item[$({\rm I}_{5.1})$] $R=R_1+R_2+R_3$, where $R_1$, $R_2$, and $R_3$ are germs of smooth curves $(M_1,R_1)_{o'}=(M_1,R_2)_{o'}=(M_1,R_3)_{o'}=1$ and the branch curve $B=F(R)$ consists of three irreducible germs,
\item[$({\rm I}_{5.2})$] $R=R_1+R_2+R_3$, where $R_1$, $R_2$, and $R_3$ are germs of smooth curves $(M_1,R_1)_{o'}=(M_1,R_2)_{o'}=(M_1,R_3)_{o'}=1$, the branch curve $B=F(R)$ consists of two irreducible germs and has singularity of type $A_{2n-1}$ for some $n\in\mathbb N$.
\end{itemize}

Note that the case, when  $R=R_1+R_2+R_3$, where $R_1$, $R_2$, and $R_3$ are germs of smooth curves ,  $(M_1,R_1)_{o'}=(M_1,R_2)_{o'}=(M_1,R_3)_{o'}=1$, and the germ of braanch curve $B=F(R)$ is irreducible, is impossible, since
$\deg_{o'} F=4$.
\begin{prop} \label{I_1.} In case $({\rm I}_1)$,  the germ of a finite cover $F$ is equivalent to the germ $F_{{\bf 1}_1}$.
\end{prop}
\proof A germ $F$ is ramified along $R_1$ with multiplicity $4$ and
the function $u=z$ is a local parameter on $R_1$ at the point $o'$. Therefore $\deg F_{\mid R_1}=1$ and $B=F(R_1)$ is a smooth curve germ  at the point $o$. Hence,  $\pi_1^{loc}(B,o)\simeq \mathbb Z$. Since the germ $F$ is ramified along $R_1$ with multiplicity $4$, then the monodromy group $G_F$ of $F$ is  a cyclic group of the fourth order and $G_F$ is generated in $\mathbb S_4$ by a cycle of length $4$. But all such subgroups are conjugated in $\mathbb S_4$. Therefore, by  Grauert - Remmert - Riemann - Stein Theorem, the germ $F$ is equivalent to the germ of finite cover given by functions $ u=z$, $v=w^4$. \qed

\begin{prop} \label{I_2} In case $({\rm I}_2)$, the germ of a finite cover $F$ is equivalent to the germ $F_{{\bf 4}_1,n}$ for some $n\in \mathbb N$.
The singularity type of the germ $(B,o)$ of the branch curve of $F_{{\bf 4}_1,n}$ is $A_{8n-1}$ and the monodromy group $G_{F_{{\bf 4}_1,n}}=\mathbb S_4$.
\end{prop}
\proof
The germ of a finite cover $F$ is ramified along $R_1$ with multiplicity $3$ and along $R_2$ with multiplicity $2$. The function $u=z$ is a local parameter on both $R_1$ and $R_2$ at the point $o'$. Therefore $B_1=F(R_1)$ and $B_2=F(R_2)$ are smooth curve germs  at the point $o$ and $\deg F_{\mid R_1}=\deg F_{\mid R_2}=1$. Note that $B_1\neq B_2$, since $\deg_{o'} F=4$. Let $(B_1,B_2)_{o}=k$. Therefore
the point $o$ is the singular point of the branch curve $B=B_1\cap B_2$ of type $A_{2k-1}$. By Lemma \ref{cl1.2},
$\pi_1^{loc}(B,o)=\langle \gamma_1, \gamma_2\, \, \mid \,\, (\gamma_1\gamma_2)^k=(\gamma_2\gamma_1)^k\rangle, $
where $\gamma_i$ are bypasses around $B_i$, $i=1,2$, and
the monodromy group $G_F$ is generated by a cycle $\tau_1=F_*(\gamma_1)$ of length $3$ and a transposition $\tau_2=F_*(\gamma_2)$. Since $G_F$ is a transitive subgroup of $\mathbb S_4$,  then $G_F=\mathbb S_4$ and with up to renumbering, we can assume that
$\tau_1=(1,2,3)$ and $\tau_2=(3,4)$. We have $\tau_1\tau_2=(1,2,4,3)$ and $\tau_2\tau_1=(1,2,3,4)$. Therefore
$(\tau_1\tau_2)^k=(\tau_2\tau_1)^k$
if and only if $k=4n$ and, by  Grauert - Remmert - Riemann - Stein Theorem, the finite cover $F$, with up to  changes of coordinates, coincides with the germ of finite cover $F_{{\bf 4}_1,n}$ given by
$$ u=z,\qquad v= w^4-z^nw^3.$$
Indeed,  the ramification divisor $R$ of $F_{{\bf 4}_1,n}$ is given by equation $4w^3-3z^nw^2=0$ and hence, $R_1$ is given by $w=0$ and $R_2$ is given by equation $4w=3z^n$. Then it is easy to see that $B_1$ is given by equation $v=0$ and $B_2$ is given by equation $4^4v+3^3u^{4n}=0$, i.e., the point $o$ is the singular point of the branch curve $B$ of type $A_{8n-1}$. To complete the proof, note that, by Proposition \ref{Anrig}, the singularities of type $A_m$ are rigid. \qed

\begin{prop} \label{I_3} In case $({\rm I}_3)$, the germ $F$ is deformation equivalent to the germ $F_{{\bf 4}_2,n,\beta}$ for some $n\in \mathbb Z_{\geq 0}$ and $\beta\neq 0$.

The germs $F$ deformation equivalent to the germ $F_{{\bf 4}_2,0,1}$ are equivalent.
\end{prop}
\proof We need the following Lemma.
\begin{lem}\label{T3} Let an irreducible curve germ $R_1\subset (U,o')$ be given by equation
\begin{equation}\label{I3} w_0^3(1+\widetilde h_3(z_0))+\sum_{i=0}^2w_0^iz_0^{k_i}(b_i+\widetilde h_i(z_0))
+w_0^4\widetilde h_4(z_0,w_0)=0,\end{equation}
where  $\widetilde h_4(z_0,w_0)\in \mathbb C[[z_0,w_0]]$, $\widetilde h_i(z_0)\in\mathfrak{m}\cap\mathbb C[[z_0]]$ for $i=0,1,2,3$, and $k_i\in \mathbb N$, $b_i\in \mathbb C^*$ for $i=0,1,2$. Then there is a local change of coordinates in $(U,o')$ of the form
$$ z=z_0,\qquad w=w_0+g(z_0),$$
where $g(z_0)\in \mathbb C[[z_0]]$ such that the equation of $R_1$ has the following form
\begin{equation}\label{I3-1} \begin{array}{c}w^3(1+h_3(z))+w^2z^{n_2}(a_2+ h_2(z))+
wz^{2n_1+\beta_1}(a_1+ h_1(z))+ \\ z^{3n_0+\beta_0}(a_0+ h_0(z))
+  w^4h_4(z,w)=0,\end{array}\end{equation}
where   $h_4(z,w)\in \mathbb C[[z,w]]$, $h_i(z)\in\mathfrak{m}\cap\mathbb C[[z]]$ for $i=0,1,2,3$, and $n_i\in \mathbb Z_{\geq 0}$, $a_i\in \mathbb C^*$ for $i=0,1,2$, and
\begin{equation} \label{cond} \beta_0=1\,\, \text{or}\,\, 2,\,\,\, \beta_1=0\,\, \text{or}\,\, 1,\,\,\,n_1\geq n_0\,\, \text{and}\,\,  \beta_1>0\,\,\,  \text{if}\,\, n_1=n_0,\,\,   n_2\geq n_0+1.\end{equation}
\end{lem}
\begin{rem} \label{remark1} Note that if we change the coordinate $w_0$ by $w_1=w_0+g(z_0)$, then to calculate the Jacobian $J_F$ of the germ of finite cover $F$ in coordinates $(z_0, w_1)$, it is sufficient to substitute $w_1-g(z_0)$  in the equation {\rm (\ref{I3})} instead of $w_0$.
\end{rem}

\noindent {\it Proof} of Lemma \ref{T3}. Let $k_0= 3m_0+\beta_0$, where $\beta_0$ is the remainder of $k_0$ divided by $3$; $k_1= 2m_1+\beta_1$, where $\beta_1$ is the remainder of $k_1$ divided by $2$; and $k_2:=m_2$.

Let us show that $m_0\leq \min(m_1,m_2)$, since otherwise, the curve germ $R_1$ is not irreducible. Indeed, let $m_1<m_0$ and
$m_1\leq m_2$  (the case when $m_2<m_0$ and $m_2\leq m_1$ is similar and therefore its consideration will be omitted). Consider a sequence of $\sigma$-processes $\sigma: U_{m_1}\to U$ given by functions $z_0=z_{m_1}$ and $w_0=w_{m_1}z_{m_1}^{m_1}$. Then the exceptional curve $E_{m_1}$ of the last $\sigma$-process is given by equation $z_{m_1}=0$ and the proper inverse
image $\sigma^{-1}(R_1)$ of $R_1$ is given by equation
$$w_{m_1}^3(1+\widetilde h_3(z_{m_1}))+\sum_{i=0}^2w_{m_1}^iz_{m_1}^{k_i-(3-i)m_1}(b_i+\widetilde h_i(z_{m_1}))
+w_{m_1}^4z_{m_1}^{m_1}\widetilde h_4(z_{m_1},w_{m_1}z_{m_1}^{m_1})=0.$$
Rewrite this equation in the following form:
\begin{equation} \label{I31} w_{m_1}^3 +b_0z_{m_1}^{3(m_0-m_1)+\beta_0}+b_1w_{m_1}z_{m_1}^{\beta_1}+b_2w_{m_1}^2z_{m_1}^{m_2-m_1}+\,\,
\text{monomials of degree}\,\,\geq 4=0.\end{equation}
It follows from (\ref{I31}) that if $\beta_1=m_2-m_1=0$ then $\sigma^{-1}(R_1)$ intersects with $E_{_1}$ at three points
$(z_{m_1},w_{m_1})=(0,0)$ and $(z_{m_1},w_{m_1})=(0,\frac{-b_1\pm \sqrt{b_1^2-4b_2}}{2})$; if $\beta_1=1$ and $m_2-m_1=0$ then $\sigma^{-1}(R_1)$ intersects with $E_{m_1}$ at two points $(z_{m_1},w_{m_1})=(0,0)$ and $(z_{m_1},w_{m_1})=(0,-b_2)$; if $\beta_1=0$ and $m_2-m_1> 0$ then $\sigma^{-1}(R_1)$ intersects with $E_{m_1}$ at three points
$(z_{m_1},w_{m_1})=(0,0)$ and $(z_{m_1},w_{m_1})=(0,\pm\sqrt{-b_1})$, which contradicts irreducibility of the curve germ $R_1$.
Finally, if $\beta_1=1$ and $m_2-m_1> 0$, then the quadratic part of the left side of  equation (\ref{I31}) is $b_1z_{m_1}w_{m_1}$, which also contradicts irreducibility of the curve germ $R_1$.

Now, let us show that we can assume that $\beta_0\neq 0$. Indeed, suppose that $\beta_0=0$.  As above, consider a sequence of $\sigma$-process $\sigma: U_{m_0}\to U$ given by functions $z_0=z_{m_0}$ and $w_0=w_{m_0}z_{m_0}^{m_0}$. Then the exceptional curve $E_{m_0}$ of the last $\sigma$-process is given by equation $z_{m_0}=0$ and the proper inverse
image $\sigma^{-1}(R_1)$ of $R_1$ is given by equation
$$w_{m_0}^3 +b_0+b_1w_{m_0}z_{m_0}^{2(m_1-m_0)+\beta_1}+b_2w_{m_0}^2z_{m_0}^{m_2-m_0}+\,\,
\text{monomials of degree}\,\,\geq 4=0.$$ 
It follows from irreducibility of the curve germ $\sigma^{-1}(R_1)$ that
$$w_{m_0}^3 +b_0+b_1w_{m_0}z_{m_0}^{2(m_1-m_0)+\beta_1}+b_2w_{m_0}^2z_{m_0}^{m_2-m_0}=(w_{m_0}+\sqrt[3]{b_0})^3.$$
Therefore, if we change the coordinate $w_0$ by $\overline w_0=w_0+\sqrt[3]{b_0}z_0^{m_0}$, then we obtain that $R_1$ is given by an equation of the same type as the equation (\ref{I3}):
$$ \overline w_0^3(1+\overline h_3(z_0))+
\sum_{i=0}^2\overline w_0^iz_0^{\overline k_i}(\overline b_i+\overline h_i(z_0))
+\overline w_0^4\overline h_4(z_0,\overline w_0)=0$$
in which $\overline k_0>k_0$.
After a finite number (since the singular point $o'$ of $R_1$ can be resolved by finite number of $\sigma$-processes) of such coordinate substitutions, we obtain that the germ $R_1$ is given by an equation of the same type as equation (\ref{I3-1})
in which $\beta_0=1$ or $2$ and $n_0\leq \min(n_1,n_2)$.

Note that the same arguments as above (the irreducibility of $R_1$ and possibility to carry out $n_1$ or $n_2$ $\sigma$-processes if necessary) show that $n_2\geq n_0+1$ and $\beta_1>0$ if $n_1=n_0$. \qed \\

In the case $({\rm I}_3)$, the ramification divisor $R=R_1$ of the germ of a finite cover $F$ is given by equation (\ref{I3}) and by Lemma \ref{T3}, we can assume that it is given by equation (\ref{I3-1}) satisfying conditions (\ref{cond}).

Then, by Remark \ref{remark1}, the germ $F$ is given by functions
\begin{equation} \label{fam1} \begin{array}{rrl} u & = & z, \\ v & = & w^4(1+h_3(z))+4wz^{3n_0+\beta_0}(a_0+ h_0(z))+ \\
& & \frac{4}{3}w^3z^{n_2}(a_2+ h_2(z))+
2w^2z^{2n_1+\beta_1}(a_1+ h_1(z)) +  w^5g(z,w),\end{array}\end{equation}
where $w^5g(z,w)= \int{w^4h_4(z,w)dw}$.

To find the singularity  type of the branch locus of $F$, let us show that $R_{1}$ can be given by parametrisation
\begin{equation} \label{D11}
z  =  \tau^3, \quad
w  =  \tau^{3n_0+\beta_0}(\sqrt[3]{-a_0}+\tau g_1(\tau)), 
\end{equation}
where $g_1(\tau)\in \mathbb C[[\tau]]$ .
Indeed, denote by $\nu: \widetilde R_{1}\to R_{1}$ the  resolution of singularity of the germ $R_{1}$. It follows from equation (\ref{I3-1}) 
that $(M_1,R_{1})_{o'}=3$ and $(M_2,R_{1})_{o'}=3n_0+\beta_0$. Therefore there exist a local parameter $\tau$ at the point $\nu^{-1}(o')$ in $\widetilde R_{1}$ such that $\nu^{-1}(z)=\tau^3$ and $\nu^*(w)=\tau^{3n_0+\beta_0}
\sum_{i=0}^{\infty} c_i\tau^i$.
If we substitute $\tau^3$ instead of $z$ and ${3n_0+\beta_0}\sum_{i=0}^{\infty}c_i\tau^{i+\beta_0}$ instead of $w$ in (\ref{I3-1}) then, as a result, we must obtain a power series identically equal to zero. In particular, we obtain that $c_0=\sqrt[3]{-a_0}$.

By (\ref{cond}), we have $n_2-n_0>0$ and $n_1-n_0+\beta_1>0$. Therefore, if we substitute  $\tau^{3n_0+\beta_0}(\sqrt[3]{-a_0}+\tau g_1(\tau))$ instead of $w$ and $\tau^3$  instead of $z$ in (\ref{fam1}), then we obtain that $B=F(R_{1})$ is given parametrically by equations of the form
\begin{equation} \label{equi1} u=\tau^3,\qquad v=\tau^{12n_0+4\beta_0}(3a_0\sqrt[3]{-a_0}+ \tau g_2(\tau)),
\end{equation}
where $g_2(\tau)\in \mathbb C[[\tau ]]$. Therefore the branch curve $B$ of $F$ has the singularity at $o$ of type $T_{3}(4n_0+\beta_0,\beta_0)$.

It is easy to see that $\deg F_{\mid R_1}=1$ and $F$ is ramified along $R_1$ with multiplicity two, therefore the monodromy group $G_F$ of the germ $F$ is generated by transpositions, and since $G_F$ is a transitive subgroup of $\mathbb S_4$, we have $G_F=\mathbb S_4$.

To complete the proof of Proposition \ref{I_3}, it suffices to consider the germs $F_{{\bf 4}_2,n_0,\beta_0}$ given by functions
$ u  =  z$, $v  =  w^4+4wz^{3n_0+\beta_0}$ and apply Lemma \ref{Lemx} and Corollaries \ref{eqi-cov1} and \ref{eqi-cov2}.  \qed \\

The following two Propositions will be proved simultaneously.
\begin{prop} \label{I_4.1} In case $({\rm I}_{4.1})$, the germ of finite cover $F$ is deformation equivalent either to the germ $F_{{\bf 4}_3,k,m}$, or to the germ $F_{{\bf 4}_4,k,m}$ for some $m\in \mathbb N$ and  $k\in \mathbb Z_{\geq 0}$.

The  singularity type of $(B_{{\bf 4}_3,k,m},o)$ is $T_{3,8m+6k+2, 8m}$ and the singularity type of $(B_{{\bf 4}_4,k,m},o)$ is $T_{3,8m+2k+4, 8m+4}$; in both cases $G_F=\mathbb S_4$.
\end{prop}

\begin{prop} \label{I_4.2} In case $({\rm I}_{4.2})$, the germ $F$ is equivalent to the germ $F_{{\bf 3}_1,2k+1}$ for some
$k\in \mathbb Z_{\geq 0}$. The singularity type of the germ $(B_{{\bf 3}_1,2k+1},o)$
is $A_{8k+3}$ and $G_{F_{{\bf 3}_1,2k+1}}$ is a dihedral group $\mathbb D_4\subset\mathbb S_4$.\end{prop}
\proof
By Remark \ref{remark0}, we can choose  coordinates $(z,w)$ in $(U,o')$ such that $w=0$ is an equation of $R_2$ and  $u=z$, where $(u,v)$ are coordinates in $(V,o)$. Since the germ $R_1$ is irreducible and $(M_1,R_1)_{o'}=2$, then arguments similar to that which we used in the proof of Lemma \ref{T3}, imply  that an equation of $R_1$ has the following form
\begin{equation}\label{eqR2}
(w-a_0z^m(1+zh_0))^2-z^{2(m+n)+1}(1+zh_1)+wz^{m+n+1}h_2+w^2zh_3+\sum_{i=3}^{\infty}w^ih_{i+1}=0, \end{equation}
where $h_i\in\mathbb C[[z]]$ for $i\in \mathbb N$, $h_0\in \mathbb C[z]$, $\deg h_0\leq n-1$,  $a_0\in \mathbb C$, $m\in\mathbb N$ and
$n\in \mathbb Z_{\geq 0}$. Therefore
$$\begin{array}{ll} \displaystyle \frac{\partial v}{\partial w}= & 4[w^3-2a_0z^m(1+zh_0)w^2+a_0^2z^{2m}(1+zh_0)^2w-z^{2(m+n)+1}w - \\
& -z^{2(m+n)+2}h_1w-z^{m+n+1}h_2w^2+zh_3w^3+ \displaystyle \sum_{i=3}^{\infty}h_{i+1}w^{i+1}] \end{array}$$
and hence
\begin{equation}\label{vI4} \begin{array}{ll} \displaystyle v = &  w^4-\frac{8a_0}{3}z^m(1+z_0h_0)w^3+2a_0^2(z^{2m}(1+zh_0)^2-4z^{2(m+n)+1})w^2 - \\
& -2z^{2(m+n)+2}h_1w^2+\frac{4}{3}z^{m+n+1}h_2w^3+zh_3w^4+
{\displaystyle \sum_{i=3}^{\infty}}\frac{4}{i+2}h_{i+1}w^{i+2}. \end{array}\end{equation}

The germ of the ramification curve $R=R_1\cup R_2$ has the singularity type $T_{3,2m+2n,2m}$ if $a_0\neq 0$ and $T_{3,2m+2n,2m+2n+1}$ if $a_0=0$.

It follows from (\ref{vI4}) that $B_2=F(R_2)$ is given by equation $v=0$.

Consider the case when $a_0\neq 0$. As in the proof of Proposition \ref{I_3}, it is easy to show that equation (\ref{eqR2}) implies that $R_1$ can be given parametrically by functions of the following form:
\begin{equation} \label{parR2}
\begin{array}{ll} z = & \tau^2, \\
w= & a_0\tau^{2m}(1+\tau^2h_0(\tau^2))+\tau^{2m+2n+1} + \displaystyle \sum_{i=2(m+n+1)}^{\infty}c_i\tau^i \end{array}\end{equation}

To obtain a parametrisation of $B_1$, we substitute $z(\tau)$ and $w(\tau)$ from (\ref{parR2}) instead of $z$ and $w$ in (\ref{vI4}) and obtain that $B_1$ has a parametrisation of the following form
\begin{equation}\label{pqrB2} \begin{array}{ll}
u= & \tau^2, \\
v= & \displaystyle  \frac{a_0^4}{3}\tau^{8m}+\displaystyle \tau^{8m}\sum_{i=1}^{\infty}a_i\tau^i
\end{array}\end{equation}
If $a_i=0$ for all odd $i$, then $\deg F_{\mid R_1}=2$, $F(R_1)$ is a smooth germ, and $B$ has the singularity type $A_{8m-1}$. Therefore $F_*(\gamma_1)\subset \mathbb S_4$ is a product of two commuting transpositions and $F_*(\gamma_2)$ is a transposition, where $\gamma_1,\gamma_2$ is a good geometric base in which $\gamma_1$ is a bypass around $B_1$ and $\gamma_2$ is a bypass around $B_2$. Then, by Lemma \ref{AnD4}, the monodromy group $G_F$ is a dihedral group $\mathbb D_4\subset \mathbb S_4$. But, Remark \ref{remx} (see the end of Subsection \ref{secti}) claims that it is impossible, since the germ $R$ consists of two irreducible germs.

Let $i_0$ be the smallest number for which $a_{i_0}\neq 0$. Then $\deg F_{\mid R_1}=1$ and $B$ has the singularity of type $T_{3,8m+i_0-1,8m}$. It follows from Lemma \ref{cl8} that $i_0=6k+3$ for some $k\in \mathbb Z_{\geq 0}$.

Let us show that in this case the germ $F$ is deformation equivalent to the germ $F_{{\bf 4}_3,k,m}$ given by functions
\begin{equation}\label{vI41} \begin{array}{ll} u= & z \\ \displaystyle v = & w^4-\frac{8}{3}z^mw^3+2a_0^2(z^{2m}
-z^{2(m+k)+1})w^2. \end{array}\end{equation}
Indeed, we have
$$\begin{array}{ll} \displaystyle \frac{\partial v}{\partial w}= & 4[w^3-2z^mw^2+z^{2m}w-z^{2(m+k)+1}w]\end{array}$$
and hence, $R_{m,k}=R_1\cup R_{2}$,  where an equation of $R_2$ is $w=0$ and an equation of $R_{1}$ is
$$(w-z^m)^2-z^{2(m+k)+1}=0.$$
The germ $R_1$ has the following parametrisation
\begin{equation}\label{parR1t}
z=\tau^2,\qquad w=\tau^{2m}(1-\tau^{2k+1}) \end{equation}
and hence,  the curve $R_{m,k}$ has a singularity of type $T_{3,2m+2k,2m}$.
The branch curve $B_{m,k}=B_1\cup B_{2}$, where $B_2$ is given by equation $v=0$. In order to obtain a parametrisation of
$B_{1}$, we substitute $z(\tau)$ and $w(\tau)$ from (\ref{parR1t}) instead of $z$ and $w$ in (\ref{vI41}). As a result, we obtain that $B_1$ is given parametrically by functions
\begin{equation}\label{pqrB2t} \begin{array}{ll}
u= & \tau^2, \\
v= & \displaystyle  \tau^{8m}(1+\tau^{2k+1})^2[(1+\tau^{2k+1})^2-\frac{8}{3}(1+\tau^{2k+1})+2-2\tau^{4k+2}]=\\
& \frac{1}{3}\tau^{8m}(1-2\tau^{4k+2}-8\tau^{6k+3}-6\tau^{8k+4}).
\end{array}\end{equation}
Consequently, the curve germs $B_{m,k}$ has the singularity of type $T_{3,8m+6k+2,8m}$. It follows from Lemma \ref{cl8} and Corollary \ref{eqi-cov1} that  $F$ is deformation equivalent to $F_{{\bf 4}_3,k,m}$.

Now, consider the case $a_0=0$ (in this case we  put $k:=m+n$). Then $R=R_1\cup R_2$ has the singularity of type
$T_{3, 2k, 2k+1}$. Let us write equation (\ref{eqR2}) of $R_1$ in the following form
\begin{equation}\label{equR2}
w^2-z^{2k+1} +\sum_{m=4k+3}^{\infty}(\sum_{2i+(2k+1)j=m} a_{i,j}z^iw^j)=0. \end{equation}
Then it is easy to see that a parametrisation of $R_1$ has the following form
\begin{equation} \label{parameR2}\begin{array}{rl}
z= & \tau^2, \\
w= & \tau^{2k+1} +\sum_{n=2k+2}^{\infty}r_n\tau^n, \end{array}. \end{equation}

We have
$$\frac{\partial v}{\partial w}=4[w^3-z^{2k+1}w +\sum_{m=4k+3}^{\infty}(\sum_{2i+(2k+1)j=m} a_{i,j}z^iw^{j+1})]$$
and consequently, $F$ is given by functions
\begin{equation} \label{F4fin} \begin{array}{ll} u= & z,\\
v=& w^4-2z^{2k+1}w^2 +\displaystyle w^2\sum_{m=4k+3}^{\infty}(\sum_{2i+(2k+1)j=m}\frac{4}{j+2} a_{i,j}z^iw^{j}).
\end{array}\end{equation}

To obtain a parametrisation of $B_1$, we substitute $z(\tau)$ and $w(\tau)$ from (\ref{parameR2}) instead of $z$ and $w$ in (\ref{F4fin}) and obtain that $B_1$ has a parametrisation of the following form
\begin{equation} \label{B2fin} \begin{array}{ll} u= & \tau^2,\\
v= & -\tau^{8k+4} +\displaystyle \sum_{m=8k+5}^{\infty}b_m\tau^m.
\end{array}\end{equation}

There are two possibilities: either there is $m_0=8k+2(n+2)+1$ being the smallest odd number such that $b_{m_0}\neq 0$ or  $b_m=0$  for all odd $m\geq 8k+5$.

Let $b_{8k+2(n+2)+1}\neq 0$. Then $B=B_1\cup B_2$ has the singularity of type $T_{3,8k+2n+4,8k+4}$, $k$
and $n\in \mathbb Z_{\geq 0}$.
Let us show that in this case the germ $F$ is deformation equivalent to the germ $F_{{\bf 4}_4,n,k}$ given by functions
\begin{equation}\label{Fkn}\begin{array}{ll}
u= & z, \\
v=& w^4-2z^{2k+1}w^2 +12z^{k+n+1}w^3.
\end{array}\end{equation}
The ramification divisor of $F_{{\bf 4}_4,n,k}$ is $R=R_1\cup R_2$, where $R_2$ is given by equation $w=0$ and $R_1$ is given by equation
\begin{equation}\label{exR2}
w^2-z^{2k+1}+3z^{k+n+1}w=0.
\end{equation}
It is easy to check that $R_1$ is given parametrically by functions
\begin{equation} \label{x} \begin{array}{ll}
z=& \tau^2, \\
w=&\tau^{2k+1} -\frac{3}{2}\tau^{2(k+n+1)} +\tau^{2(k+n)+3}h(\tau)),
\end{array}\end{equation}
where $h(\tau)\in\mathbb C[[\tau]]$.
To obtain a parametrisation of $B_1=F_{{\bf 4}_4,n,k}(R_1)$, let us substitute $z(\tau)$ and $w(\tau)$ from (\ref{x}) in (\ref{Fkn}) and, as a result, we obtain \begin{equation} \label{xB} \begin{array}{ll}
u=& \tau^2, \\
v=&-\tau^{8k+4} +6\tau^{8k+2n+5} +\tau^{8k+2n+6}h_1(\tau),
\end{array}\end{equation}
where $h_1(\tau)\in\mathbb C[[\tau]]$. Consequently, $B$ has the singularity type $T_{3,8k+2n+4, 8k+4}$ and Proposition \ref{I_4.1} follows from Lemma \ref{cl7} and Corollary \ref{eqi-cov1}.

Consider the case when $b_m=0$  for all odd $m\geq 8k+5$. It follows from (\ref{B2fin}) that $\deg F_{\mid R_1}=2$,
$B_1=F(R_1)$ is a smooth germ, and the singularity type of the germ $(B,o)$ is $A_{8k+3}$. Therefore $F_*:\pi_1^{loc}(B,o)\to\mathbb S_4$ has the following properties: $(*_{1})$ $G_F=\text{im}, F_*$ is a transitive group of $\mathbb S_4$;
$(*_{2})$ $F_*(\gamma_1)$ is the product of two commuting transpositions and $F_*(\gamma_2)$
is a transposition, where $\gamma_1$ is a bypass around $B_1$ and $\gamma_2$ is a bypass around $B_2$.

By Lemma \ref{AnD4}, $G_F=\mathbb D_4$ is a dihedral subgroup of $\mathbb S_4$.
Therefore, 
Corollary \ref{eqi-cov1}, to complete the proof of Proposition \ref{I_4.2}, it suffices to show that the branch curve $B$ of the germ of finite cover $F_{{\bf 3}_1,2k+1}$ given by functions
\begin{equation} \label{I42x}
u=z, \quad v=w^4-2z^{2k+1}w^2,
\end{equation}
has the singularity of type $A_{8k+3}$ and $F_{{\bf 3}_1,2k+1*}$ has properties $(*_{1})$ and $(*_{2})$.

The ramification divisor $R=R_1\cup R_2$ of $F_{{\bf 3}_1,2k+1}$ is given by 
$w(w^2-z^{2k+1})=0$, where $R_1$ is given by $w^2-z^{2k+1}=0$.
Therefore  $R_1$ has the following parametrisation:
$$ u=\tau^2,\qquad v=\tau^{2k+1}.
$$
Consequently, $F_{{\bf 3}_1,2k+1\mid R_1}:R_1\to (V,o)$  is given by functions
\begin{equation}\label{B1111} u=\tau^2,\qquad v=-\tau^{8k+4}\quad (=(\tau^{2k+1})^4-2(\tau^2)^{2k+1}(\tau^{2k+1})^2)
\end{equation}
and, hence $\deg F_{{\bf 3}_1,2k+1\mid R_1}=2$, $B_1$ is a smooth germ touching the germ $B_2$ given by $v=0$ with multiplicity $4k+2$. Therefore
the branch curve $B$ of the germ $F_{{\bf 3}_1,2k+1}$ has the singularity of type $A_{8k+3}$ and $F_{{\bf 3}_1,2k+1*}$ has properties $(*_1)$ and $(*_2)$.\qed \\

The following two Propositions also will be proved simultaneously.

\begin{prop} \label{I5.1} In case $({\rm I}_{5.1})$, the germ $F$ is deformation equivalent either to one of the germs $F_{{\bf 4}_5,k,m}$, $ m,k\geq 1$, {\rm (}the singularity type  of $B_{F_{{\bf 4}_5,k,m}}$ is $T_{3,8m+6k-1,8m}${\rm )}, or
to one of the germs $F_{{\bf 4}_6,k,m}$, $m,k\geq 1$, {\rm (}the singularity type of $B_{F_{{\bf 4}_6,k,m}}$
is $T_{3,8m+2k-1,8m}${\rm )}, or to one of the germs $F_{{\bf 4}_7,m}$, $m\geq 1$, {\rm (}the singularity type of  $B_{F_{{\bf 4}_7,m}}$ is $T_{3,8m-1,8m}${\rm )}. In all cases, the monodromy groups of the germs of finite covers are $\mathbb S_4$.
\end{prop}

\begin{prop} \label{I5.2}
In case $({\rm I}_{5.2})$, the germ $F$ is equivalent to the germ $F_{{\bf 3}_1,2k}$ for some $k\in \mathbb N$.
The singularity type of the germ of branch curve $(B_{F_{{\bf 3}_1,2k}},o)$ is $A_{8k-1}$ and the monodromy group $G_{F_{{\bf 3}_1,2k}}$ is the dihedral group $\mathbb D_4\subset\mathbb S_4$.\end{prop}

\proof In cases $({\rm I}_{5.1})$ and $({\rm I}_{5.2})$ the ramification divisor $R$  consists of three irreducible germs,
$R=R_1\cup R_2\cup R_3$. Let us renumber them so that
$$(R_1,R_3)_{o'}=(R_2,R_3)_{o'}=m\leq (R_1,R_2)_{o'}=m+n$$
and choose coordinates $z,w$ such that $w=0$ is an equation of $R_3$.

Since $(M_1,R_i)_{o'}=1$, then $z$ is a local parameter on each germ $R_i$ and we can choose equations of $R_i$ of the following form:  $w=f_i(z)$, $i=1,2$, where
if $n\geq 1$ then
\begin{equation}\label{equaI5y}
f_1(z)=z^mp(z)+ z^{m+n}h_1(z), \qquad f_2(z)=z^mp(z)+ z^{m+n+n_2}h_2(z),
\end{equation}
where $n_2\geq 0$, and if $n=0$ then
\begin{equation}\label{equalI5z} f_i(z)=z^mh_i(z),\end{equation}
where $p(z)\in \mathbb C[z]$, $\deg p(z)\leq n-1$, $p(0)\neq 0$ and $h_i(z)\in \mathbb C[[z]]$, $ h_1(0)\neq 0$, $h_2(0)\neq 0$, and $h_1(0)\neq h_2(0)$ if $n_2=0$.
If $n\geq 1$ then
$$\frac{\partial v}{\partial w}=4w(w-z^mp(z)- z^{m+n}h_1(z))(w-z^mp(z)- z^{m+n+n_2}h_2(z))(1+g_1),$$
where $g_1\in \mathfrak{m}$. Consequently, $F$ is given by functions
\begin{equation}\label{part}\begin{array}{l} u=  z \\ v=  w^4(1+g_4)-\frac{4}{3}[2z^mp(z)+z^{m+n}h_1(z)+z^{m+n+n_2}h_2(z))](1+g_3)w^3+\\ 2[p(z)^2z^{2m}+p(z)(h_1(z)z^{2m+n}+h_2(z)z^{2m+n+n_2})+h_1(z)h_2(z)z^{2(m+n)+n_2}](1+g_2)w^2,\end{array}
\end{equation}
where $g_i\in\mathfrak{m}$. Similarly, if $n=0$ then $F$ is given by functions
\begin{equation}\label{part1}\begin{array}{l} u=  z \\ v=  w^4(1+g_4)-\frac{4}{3}z^{m}[h_1(z)+h_2(z)](1+g_3)w^3+ 2z^{2m}h_1(z)h_2(z)(1+g_2)w^2.\end{array}
\end{equation}

We have $v=0$ is an equation of $B_3=F(R_3)$. To obtain equations of $B_i=F(R_i)$, let us substitute $f_i(u)$ from (\ref{equaI5y}) or (\ref{equalI5z}), resp., in (\ref{part}) or (\ref{part1}) instead of $w$. After substitutions, we obtain that
the equation of $B_i$, $i=,1,2$, have the following form
\begin{equation}\label{eqBi1}
v=\frac{p(u)^4}{3}u^{4m}(1 + uh_{i,3}(u))
\end{equation}
if $n>0$ and
\begin{equation}\label{eqBi2}
v=\frac{h_i(u)^3}{3}(2h_j(u)-h_i(u))u^{4m}(1+uh_4(u))
\end{equation}
if $n=0$, where $\{ j\}=\{ 1,2\}\setminus \{ i\}$ and $h_{i,3}(u),h_4(u)\in \mathbb C[[u]]$.

Consider the case when $B=B_1\cup B_2\cup B_3$ consists of three different germs. Then in case $n>0$, it follows from (\ref{eqBi1}) that $(B_1,B_3)_o=(B_2,B_3)_o=4m$, and $(B_1,B_2)_o= 4m+k>4m$ for some $k$, since $p(0)\neq 0$. Therefore the singularity type of $(B,o)$ is $T_{3,8m+2k-1,8m}$.

In particular, it is  easily to check that the ramification locus of the germ $F_{{\bf 4}_5,k,m}$ given by functions
$$u=z,\quad v=w^4-\frac{8}{3}z^m(1+2z^k)w^3+2z^{2m}(1+4z^k+3z^{2k})w^2,\qquad m,k\geq 1,$$
is given by equation $w(w-z^m-z^{m+k})(w-z^m-3z^{m+k})=0$ (i.e., $R$ has the singularity of type $T_{3,2m+2k-1,2m}$) and the singularity  type of the branch curve $B$ is $T_{3,8m+6k-1,8m}$.

Consider the case when $n=0$. Note that $2h_2(0)-h_1(0)\neq  2h_1(0)-h_2(0)$ and, in particular, it is impossible that $2h_2(0)-h_1(0)=0$ and $ 2h_1(0)-h_2(0)=0$ simultaneously. Therefore it follows from (\ref{eqBi2}) that if $2h_1(0)-h_2(0)\neq 0$ and $2h_2(0)-h_1(0)\neq 0$, then $B$ has the singularity of type $T_{3,8m-1,8m}$, and if $2h_1(0)-h_2(0)=0$ (resp., if $2h_2(0)-h_1(0)=0)$, then $B$ has the singularity type $T_{3,8m+2k-1,8m}$ for some $k\geq 1$.

In particular, it is easily to check that the ramification locus of the germ $F_{{\bf 4}_6,k,m}$ given by functions
$$u=z,\quad v=w^4-\frac{4}{3}z^m(3+z^k)w^3+4z^{2m}(1+z^k)w^2,\qquad m,k\geq 1,$$
is given by equation $w(w-z^m-z^{m+k})(w-2z^m)=0$ (i.e., $R$ has the singularity of type $T_{3,2m-1,2m}$) and the singular type of the branch curve $B$ is $T_{3,8m+2k-1,8m}$.

And again, it is easily to check that the ramification locus of the germ $F_{{\bf 4}_7,m}$ given by functions
$$u=z,\quad v=w^4-\frac{16}{3}z^mw^3+6z^{2m}w^2,\qquad m\geq 1,$$
is given by equation $w(w-z^m)(w-3z^m)=0$ (i.e., $R$ has the singularity of type $T_{3,2m-1,2m}$) and the singularity type of the branch curve $B$ is $T_{3,8m-1,8m}$.

So, in all cases considered above, the branch curves $B$ of $F$ have  singularity types $T_{3, 8m+2k-1, 8m}$ for some $k\in \mathbb Z_{\geq 0}$. The germs $F_{{\bf 4}_5,k,m}$ and $F_{{\bf 4}_6,3k,m}$ have the branch curves of the same singularity types. But, by Proposition \ref{equiD}, the covers $F_{{\bf 4}_5,k,m}$ and $F_{{\bf 4}_6,3k,m}$ are not deformation equivalent, since their ramification divisors have different types of singularity. Therefore to complete the proof of Proposition \ref{I5.1}, it suffices to apply Lemma \ref{cl181} and Corollary \ref{eqi-cov1}.

If $F(R_1)=F(R_2)$ in case when $n>0$ or if $F(R_3)$ coincides with either $F(R_1)$ or with  $F(R_2)$ in case when $n=0$, then it follows from (\ref{eqBi1}) and (\ref{eqBi2}) that $B$ consists of two irreducible germs and has the singularity type $A_{8m-1}$, since $F(R_i)$  are germs of smooth curves. Denote the irreducible components of $B$ by $B_1$ and $B_2$ and let $F^{-1}(B_2)$ be the union of two components of the divisor $R$. Then
$G_F$ is generated by a transposition $\tau_1=F_*(\gamma_1)$, where $\gamma_1$ is a bypass around $B_1$ and a product $\nu=[\tau_2\tau_1]=F_*(\gamma_2)$ of two commuting transpositions $\tau_2$ and $\tau_3$, where $\gamma_2$ is a bypass around $B_2$. Therefore, by Lemma \ref{AnD4},  $G_F$ is a dihedral subgroup $\mathbb D_4\subset \mathbb S_4$.

Let us show that in this case $F$ is equivalent to the cover $F_{{\bf 3}_1,2m}$ given by functions
$$ u=z,\qquad  v =w^4-2z^{2m}w^2.$$
Indeed, it is easy to see that the branch locus $R$ of $F_{{\bf 3}_1,2m}$ is given by equation $w(w-z^m)(w+z^m)$, the germ $B_1$ is given by equation $v=0$ and $B_2$ is given by equation $v+u^{4m}=0$. Therefore, by Lemma \ref{AnD4} and Corollary \ref{eqi-cov}, the cover $F$ is equivalent to $F_{{\bf 3}_1,2m}$. This completes the proof of Proposition \ref{I5.2}. \qed
\begin{rem} \label{remx} It follows from the proof of Proposition \ref{I5.2} that if the branch curve $B$ of the germ of finite cover $F$, $\deg_{o'} F=4$, has the singularity of type $A_{8m-1}$ and  $G_F=\mathbb D_4\subset \mathbb S_4$, then the ramification locus $R$ of $F$ consists of three irreducible components.
\end{rem}

\subsection{Case $({\rm II})$.} To comlete the proof of Theorem \ref{main}, it suffices to prove
\begin{prop} \label{II}
In case $({\rm II})$, the ramification divisor $R$ of the germ of finite cover $F$ is the union of two smooth curves, $R=R_1\cup R_2$, meeting transversally at the point $o'$ and $F$ is ramified along $R$ with multiplicity two.

There are three possibilities:
\begin{itemize}
\item[$({\rm II}_1)$] $\deg F_{\mid R_1}=\deg F_{\mid R_2}=1$, the point $o\in V$ is the singular point of types $A_{n_1}$ and $A_{n_2}$ of the curves $B_1=F(R_1)$ and $B_2=F(R_2)$ for some $n_1$, $n_2\in \mathbb N$, and $(B_1,B_2)_o=4$, the germ $F$ is deformation equivalent to the germ $F_{{\bf 4},n_1,n_2}$, the branch curve $B_{F_{{\bf 4}_8},n_1,n_2}$
 has the singularity of type $T_{4,2n_1,2n_2}$ and   $G_{F_{{\bf 4}_8, n_1,n_2}}$ is $\mathbb S_4$.
\item[$({\rm II}_2)$] accurate to the numbering, $\deg F_{\mid R_1}=1$, $\deg F_{\mid R_2}=2$, the point $o\in V$ is the singular point of the curve germ $B_1=F(R_1)$ of type $A_{2n}$ for some $n\geq 1$, the curve germ $B_2=F(R_2)$ is smooth, and $(B_1,B_2)_o=2$, the germ $F$ is
    equivalent to the germ $F_{{\bf 3}_2, n}$, the branch curve $B_{F_{{\bf 3}_2, n}}$ has the singularity of type $D_{2n+3}$, the monodromy group $G_{F_{{\bf 3}_2, n}}$ is a dihedral group
    $\mathbb D_4\subset \mathbb S_4$.
\item[$({\rm II}_3)$] $\deg F_{\mid R_1}=\deg F_{\mid R_2}=2$, the curve germs $B_1=F(R_1)$ and $B_2=F(R_2)$ are smooth, and $(B_1,B_2)_o=1$, the germ $F$ is equivalent to the germ $F_{{\bf 2}_1}$, the branch curve $B_{F_{{\bf 2}_1}}$ has the singularity of type $A_1$, the monodromy group $G_{F_{{\bf 2}_1}}$ is the Klien four group $Kl_4\subset \mathbb S_4$.
\end{itemize}
\end{prop}
\proof In case $({\rm II})$, the cover $F$ is given by functions of the following form:
$$ \begin{array}{l} u_0=f_1(z_0,w_0)f_2(z_0,w_0)+f_{\geq 3}(z_0,w_0), \\
v_0=h_1(z_0,w_0)h_2(z_0,w_0)+h_{\geq 3}(z_0,w_0),
\end{array} $$
where $f_i(z_0,w_0)$, $h_i(z_0,w_0)$, $i=1,2$, are linear forms and $f_{\geq 3}(z_0,w_0)$, $h_{\geq 3}(z_0,w_0)\in \frak{m}^3$.
Note that the forms $f_i(z_0,w_0)$ and $h_j(z_0,w_0)$ are linear independent for each pair $(i,j)$, $i,j=1,2$, since $(M_1,M_2)_{o'}=4$.

In the beginning, let us show that there are coordinates $z_1,w_1$ in $U$ and $u_1,v_1$ in $V$ such that $F$ is given by functions
\begin{equation} \label{form1} u_1 =z_1^2+g_1(z_1,w_1), \qquad v_1 =w_1^2+g_2(z_1,w_1),\end{equation}
where $g_1(z_1,w_1)$ and $g_2(z_1,w_1)\in \frak{m}^3$. Accurate to the coordinate change $u_0\leftrightarrow v_0$, we have two possibilities: 1) $f_1(z_0,w_0)$ and $f_2(z_0,w_0)$ are proportional and
2) $f_1(z_0,w_0)$ and $f_2(z_0,w_0)$ are linear independent.

In the first case we can assume (after coordinate changes in $U$ and $V$) that
$$f_1(z_0,w_0)f_2(z_0,w_0)=z_0^2\quad \text{and}\quad h_1(z_0,w_0)h_2(z_0,w_0)=w_0^2+az_0w_0$$ with some $a\in \mathbb C$. After the coordinate changes $z_0=z_1$, $w_0=\frac{1}{2}(w_1-az_1)$ and $u_0=u_1$, $v_0=4v_1$, we obtain functions of the form (\ref{form1}).

In the second case we can assume (after changes of coordinates in $U$ and $V$) that
$$f_1(z_0,w_0)=z_0,\,\, f_2(z_0,w_0)=w_0\,\, \text{and}\,\, h_1(z_0,w_0)h_2(z_0,w_0)=az_0^2+cz_0w_0+bw_0^2,$$
where $c=0$ and $a\neq 0$, $b\neq 0$
(if $c\neq 0$ then the change of coordinates  $\widetilde u_0=u_0$, $\widetilde v_0= v_0-cu_0$ "kills"\, the coefficient of $z_0w_0$). After the changes of coordinates
$$ z_1=\frac{\sqrt{a}z_0+\sqrt{b}w_0}{2},\, w_1=\frac{\sqrt{b}w_0-\sqrt{a}z_0}{2}\,\,\,  \text{and}\,\,  u=\frac{u_0+\sqrt{ab}v_0}{2\sqrt{ab}},\, v=\frac{\sqrt{ab}v_0-u_0}{2\sqrt{ab}}$$
 we obtain functions of the form (\ref{form1}).

We have
$$J(F):= \det \left(\begin{array}{cc} \frac{\partial u}{\partial z_1} & \frac{\partial u}{\partial w_1} \\ \frac{\partial v}{\partial z_1} & \frac{\partial v}{\partial w_1}\end{array}\right)= 4z_1w_1+ S(z_1,w_1), $$
where $S(z_1,w_1)\in\frak{m}^3$. Therefore
$$J(F)=4(z_1+S_1(z_1,w_1))(w_1 +S_2(z_1,w_1))$$
in the ring $\mathbb C[[z_1,w_1]]$, where $S_i(z_1,w_1)\in \frak{m}^2$ and hence $F$ is ramified along $R$ with multiplicity two.

Let us make the change of coordinates $z_2=z_1+S_1(z_1,w_1)$, $w_2=w_1+S_2(z_1,w_1)$ and write functions $u$ and $v$ in the form
$$u =z_2^2(1+z_2\widetilde f_1(z_2))+w_2\widetilde g_1(z_2,w_2), \qquad v =w_2^2(1+w_2\widetilde h_1(w_2))+z_2\widetilde g_2(z_2,,w_2),$$
where $\widetilde f_1(z_2)\in \mathbb C[[z_2]]$, $\widetilde h_1(w_2)\in \mathbb C[[w_2]]$. Note that the ramification divisor $R$ is given by  equation $z_2w_2=0$.
Therefore $R=R_1+R_2$, where $R_1=\{ z_2=0\}$ and $R_2=\{ w_2=0\}$.

Finally, if we put $z=z_2\sqrt{1+z_2\widetilde f_1(z_2)}$ and
$w=w_2\sqrt{1+w_2\widetilde h_1(w_2)}$, we obtain that $F$ is given by functions of the following form
\begin{equation} \label{form2}
u =z^2+\sum_{i=3}^{\infty}\alpha_iw^i +zw\widetilde g_3(z,w), \qquad v =w^2+\sum_{i=3}^{\infty}\beta_i z^i+zw\widetilde g_4(z,,w),
\end{equation}
where $\widetilde g_3(z,w)$ and $\widetilde g_4(z,w)\in \frak{m}$. Note that $z=0$ is an equation of $R_1$ and $w=0$ is an equation of $R_2$. Note also that $w$ is a local parameter in $R_1$ at the point $o'$ and  $z$ is a local parameter in $R_2$.

The restrictions $F_{\mid R_1}$ to $R_1$ and $F_{\mid R_2}$ to $R_2$ of $F$ are given by functions
\begin{equation} \label{R12} F_{\mid R_1}=\{ u=\sum_{i=3}^{\infty}\alpha_iw^i, \,\, v=w^2\}\,\,\, \text{and}\,\,\,  F_{\mid R_2}=\{ u=z^2,\,\,  v= \sum_{i=3}^{\infty}\beta_iz^i\} .\end{equation}
If $\alpha_i=0$ for all odd $i$, then it is easy to see that $\deg F_{\mid R_1}=2$ and $B_1=F(R_1)$ is given by equation
\begin{equation} \label{R11} u=\sum_{j=2}^{\infty} \alpha_{2j}v^j.  \end{equation}
Similarly, if $\beta_i=0$ for all odd $i$, then  $\deg F_{\mid R_2}=2$ and $B_2=F(R_2)$ is given by equation
\begin{equation} \label{R22} v=\sum_{j=2}^{\infty} \beta_{2j}u^j.  \end{equation}

If $i=2n_1+1$ is the smallest odd index for which $\alpha_i\neq 0$, then it follows from (\ref{R12}) that $\deg F_{\mid R_1}=1$, the germ  $B_1=F(R_1)$ has the singularity of type $A_{2n_1}$, and the coordinate axis $\{ u=0\}$ is the tangent line of $B_1$ at the point $o$.

Similarly, if $i=2n_2+1$ is the smallest odd index for which $\beta_i\neq 0$, then  $\deg F_{\mid R_2}=1$, the germ  $B_2=F(R_2)$ has the singularity of type $A_{2n_2}$, and the coordinate axis $\{ v=0\}$ is the tangent line of $B_2$ at the point $o$. 

If $i=2n_1+1$ and $j=2n_2+1$ are the smallest odd indexes for which $\alpha_i\neq 0$ and $\beta_j\neq 0$, then the branch curve
$B$ has the singularity of type $T_{4, 2n_1,2n_2,2}$, the monodromy group $G_F$ is $\mathbb S_4$, since $G_F$ is a transitive subgroup of $\mathbb S_4$ generated by transpositions,  and it is easy to see that the family of the germs of finite covers given by functions
\begin{equation} \label{form22}
\begin{array}{ll} u = & z^2+\alpha_{2n+1}w^{2n_1+1}+t(\displaystyle \sum_{i=1}^{\infty}\alpha_{2i} w^{2i}+\sum_{i=n_1+1}^{\infty}\alpha_{2i+1} w^{2i+1}+zw\widetilde g_3(z,,w)),
\\ v = & w^2+\beta_{2n_2+1}z^{2n_2+1}+t(\displaystyle \sum_{j=1}^{\infty}\beta_{2j} z^{2j}+\sum_{j=n_2+1}^{\infty}\beta_{2j+1} z^{2j+1}+zw\widetilde g_4(z,,w))
\end{array}\end{equation}
defines a deformation equivalence between $F$ and the gxerm $F_{{\bf 4}_8,n_1,n_2}$ given by
$$u =  z^2+\alpha_{2n_1+1}w^{2n_1+1}, \qquad v =  w^2+\beta_{2n_2+1}z^{2n_2+1}.$$

If $i=2n+1$ is the smallest odd index for which $\beta_i\neq 0$ and $\alpha_i=0$ for all odd $i$, then the branch curve
$B=B_1\cup B_2$ has the singularity of type $D_{2n+3}$  and it is easy to see that the family of the germs of finite covers given by functions
\begin{equation} \label{form21}
\begin{array}{ll} u = & z^2+\displaystyle t(\sum_{i=3}^{\infty}\alpha_iw^i +zw\widetilde g_3(z,w)), \\ v = & w^2+\beta_{2n+1}z^{2n+1}+\displaystyle t(\sum_{i=1}^{\infty}\beta_{2i} z^i+\sum_{i=n+1}^{\infty}\beta_{2i+1} z^{2i+1}+zw\widetilde g_4(z,,w))
\end{array}\end{equation}
defines a deformation equivalence between $F$ and the germ $F_{{\bf 3}_2,n}$ given by
$$u =  z^2, \qquad v =  w^2+\beta_{2n+1}z^{2n+1}.$$
The group $\pi_1^{loc}(B,o)$ is generated by bypasses $\gamma_1$ around $B_1$ and $\gamma_2$, $\gamma_3$ around $B_2$. The permutation  $F_*(\gamma_1)$ is a product of two commuting transpositions, since $\deg F_{\mid R_1}=2$, and $F_*(\gamma_i)$, $i=2,3$. are transpositions, since $\deg F_{\mid R_2}=1$. These permutations generate a transitive group in $\mathbb S_4$. It follows from Lemma \ref{Lem} that $G_{F_{{\bf 3}_2,n}}\simeq \mathbb D_4$. Recall also that the germ of singularity $D_n$ is rigid. Therefore, by Lemma  \ref{Lem} and Corollary \ref{eqi-cov}, the germ $F$ is equivalent to $F_{{\bf 3}_2,n}$.

If $\alpha_i=0$ and $\beta_i=0$ for all odd $i$, then it follows from (\ref{R11}) and (\ref{R22}) that the singular point of
$B=B_1\cup B_2$ is of type $A_1$. Consequently, $\pi_1^{loc}(B,o)=\mathbb Z\times\mathbb Z$ is generated by bypasses $\gamma_1$  and $\gamma_2$ around $B_1$ and $B_2$. The permutations $F_*(\gamma_1)$ and $F_*(\gamma_2)$ are products of two commuting transpositions, since $\deg F_{\mid R_1}=\deg F_{\mid R_2}=2$. These permutations generate a transitive group in $\mathbb S_4$. Therefore $G_F$ is the Klein four group $Kl_4\subset \mathbb S_4$. The singularity of type $A_1$ with up to change of coordinates is given by equation $uv=0$ and for the singularity $(B,o)$ of type $A_1$, there is the unique epimomorphism from  $\pi_1^{loc}(B,o)$ to $Kl_4$. Therefore the germ $F$ is equivalent to the germ $F_{{\bf 2}_1}$ given by functions $u=z^2$ and $v=w^2$.

\section{On the monofromy groups of germs of finite covers}
\subsection{On the transitive subgroups of the symmetric groups.}
Let $G$ be a finite group. We say that two subgroups $H_1$ and $H_2$ of $G$ are {\it equivalent} if there is an inner automorphism $g\in Aut(G)$ such that $g(H_1)=H_2$. A subgroup $H$ of $G$ is called {\it relatively simple} if there is not a proper non-trivial normal subgroup of $G$ contained in $H$. Denote by $A_G$ the set of representatives of equivalence classes of relatively simple subgroups of $G$ and let $I_G\subset N$ be the set of indices $i_H=(G:H)$ of subgroups $H\in A_G$.

We say that an imbedding $\varphi:G\to \mathbb S_d$ is {\it transitive} if $\varphi(G)$ is transitive subgroup of the symmetric group $\mathbb S_d$.

\begin{prop} \label{tras-subgr} The set of transitive imbedding of $G$, considered up to conjugations in symmetric groups, is in one-to-one correspondence with the set $A_G$. For any $d\in I_G$ there is a transitive imbedding $\varphi:G\to \mathbb S_d$.
\end{prop}
\proof Consider a symmetric group $\mathbb S_d$ as the group acting on the interval of natural numbers $\mathbb N_d=\{ 1,\dots, d\}$ and denote by $\mathbb S_{d-1}$ the subgroup of $\mathbb S_d$ consisting of the permutations $\tau\in\mathbb S_d$ leaving fixed $1$. Let $G$ be a transitive subgroup of $\mathbb S_d$. Then $H=G\cap \mathbb S_{d-1}$ is a subgroup of $G$ of index $(G:H)=d$, since $G$ is a transitive subgroup of $\mathbb S_d$. Let us show that $H$ is a relatively simple subgroup of $G$. Indeed, assume that a normal subgroup $N$ of $G$ is contained in $H$ and $h\in N$ is a non-trivial element. But, in this case for any $i\in \mathbb N_d$, there is an element $g_i\in G$ such that $g_i(1)=i$ and therefore $h(i)=i$ for each $i\in\mathbb N_d$, since $g_i^{-1}hg_i\in N\subset H$. As a result, we get a contradiction with the assumption that $G$ acts effectively on $\mathbb N_d$.

Conversely, let $H$ be a subgroup of $G$. Then $G$ acts on the set of left cosets of  $H$ in $G$. This action defines a homomorphism $\varphi: G\to \mathbb S_d$, where $d=(G:H)$. It is easy to check that $\varphi$ is an imbedding if and only if $H$ is a relatively simple subgroup of $G$. \qed \\

Note that for each finite group $G$ there is at least one transitive imbedding, namely, Cayley's imbedding
$c:G\hookrightarrow \mathcal S_{|G|}$ corresponding to the trivial subgroup of $G$.

\subsection{Germs of Galois smooth covers.} Consider a germ of finite cover $F:(U,o')\to (V,o)$, $\deg_{o'} F=d$, and its monodromy
$F_*:\pi_1^{loc}(B,o)\to G_F\subset \mathbb S_d$. Let $c:G=G_F\hookrightarrow \mathcal S_{|G|}$ be  Cayley's imbedding. By  Grauert - Remmert - Riemann - Stein Theorem, the homomorphism $c\circ F_*: \pi_1^{loc}(B,o)\to \mathbb S_{|G_F|}$ defines a {\it germ of Galois cover} $\widetilde F: (\widetilde U, \widetilde o)\to (V,o)$ of degree $\deg_{\widetilde o} \widetilde F=|G|$, where $\widetilde F$ is a holomorphic finite map and in general case $(\widetilde U,\widetilde o)$ is an irreducible germ of a normal complex-analytic variety. The group $G$ acts on $(\widetilde U,\widetilde o)$ such that the quotient variety
$(\widetilde U,\widetilde o)/G$ is $(V,o)$ and $\widetilde F$ is the quotient map. By Proposition \ref{tras-subgr}, the imbedding $G=G_F\hookrightarrow \mathbb S_d$ corresponds to a relatively simple subgroup $H$ of $G$, $(G:H)=d$, and it is well known that we can choose a subgroup $H_F$ of $G$ equivalent to $H$ such that the quotient variety
$(\widetilde U,\widetilde o)/H_F$ is biholomorphic to $(U,o')$ and the cover $\widetilde F$ is the composition of two covers, $\widetilde F=F\circ  F_{H_F}$, where $F_{H_F}: (\widetilde U,\widetilde o)\to (U,o')$ is the quotient map defined by the action of $H_F$ on $(\widetilde U,\widetilde o)$.

Before to formulate the statement describing the germs of the Galois smooth covers, let us recall the invariants of the binary linear groups given by Klein \cite{Kl}. For each of the three groups
(tetrahedral, octahedral, icosahedral) there are three invariant forms which we denote by $\varphi(z,w)$, $\psi(z,w)$, and $\theta(z,w)$, where $\psi(z,w)$ is the Hessian of $\varphi(z,w)$, and $\theta(z,w)$ is the Jacobian of $\varphi(z,w)$ and $\psi(z,w)$. We have
$$\begin{array}{c} \varphi_4(z,w)=z^4+2\sqrt{-3}z^2w^2+w^4, \quad \psi_4(z,w)=z^4-2\sqrt{-3}z^2w^2+w^4, \\ \theta_6(z,w)=zw(z^4-w^4)\end{array}$$
in the case of tetrahedral group;
$$\begin{array}{c}\varphi_6(z,w)=zw(z^4-w^4), \quad \psi_8(z,w)=z^8+14z^4w^4+w^8, \\ \theta_{12}(z,w)=z^{12}-33(z^8w^4+z^4w^8)+w^{12} \end{array}$$
in the case of octahedral group; and
$$\begin{array}{l}\varphi_{12}(z,w)=zw(z^{10}+11z^5w^5-w^{10}), \\ \psi_{20}(z,w)=-(z^{20}+w^{20})+288(z^{15}w^5+z^5w^{15})-494z^{10}w^{10}, \\ \theta_{30}(z,w)=(z^{30}+w^{30})+522(z^{25}w^5-z^5w^{25}) -1005(z^{20}w^{10}+z^{10}w^{20}) \end{array}$$
in the case of icosahedral group.
\begin{thm} \label{Gal} {\rm (}\cite{Sh}{\rm )} A germ of a Galois smooth cover $\widetilde F:(\widetilde U,\widetilde o)\to (V,o)$ is equivalent to one of the following rigid germs:
\begin{itemize}
\item[$1)$] $F_{1}=\{ u=z^2+zw+w^2,\,\, v=z_1^2w+zw^2\}$,
$T(B_{1})=A_2$;
\item[$2_1)$] $F_{2_1,pq,q}=\{ u=z^{pq}+w^{pq},\,\,v=z^qw^{q}\}$,
 $T(B_{2_1,pq,q})=D_{p+2}$, $p>1$, $q>1$;
\item[$2_2)$] $F_{2_2,p,1}=\{ u=z^{p}+w^p,\,\,v=zw\}$,
$T(B_{2_2,p,p})=A_{p-1}$, $p>1$, $q=1$;
\item[$2_3)$] $F_{2_3,q,q}=\{ u=z^{q}+w^{q},\,\,v=z^qw^{q}\}$,  $T(B_{2_3,q,1})=A_{3}$, $p=1$, $q>1$;
\item[$3)$] $F_{3,n_1,n_2}=\{ u=z^{n_1},\,\,v=w^{n_2}\}$, $T(B_{3,n_1,n_2})=A_1$, $n_1,n_2\in\mathbb N$;
\item[$4)$] $F_{4}=\{ u=\psi_4(z,w),\,\, v=\theta_6(z,w)\}$, $T(B_{4})=A_2$;
\item[$5)$] $F_{5}=\{ u=\varphi^3_4(z,w),\,\, v=\theta_6(z,w)\}$, $T(B_{5})=A_3$;
\item[$6)$] $F_{6}=\{ u=\varphi_4(z,w),\,\, v=\psi_4^3(z,w)\}$, $T(B_{6})=A_5$;
\item[$7)$] $F_{7}=\{ u=\varphi_4^3(z,w),\,\, v=\psi^3_4(z,w)\}$, $T(B_{7})=D_4$;
\item[$8)$] $F_{8}=\{ u=\psi_8(z,w),\,\, v=\theta_{12}(z,w)\}$, $T(B_{8})=A_2$;
\item[$9)$] $F_{9}=\{ u=\psi_8(z,w),\,\, v=\theta^2_{12}(z,w)\}$, $T(B_{9})=A_5$;
\item[$10)$] $F_{10}=\{ u=\psi^3_8(z,w),\,\, v=\theta_{12}(z,w)\}$, $T(B_{{10}})=A_3$;
\item[$11)$] $F_{11}=\{ u=\varphi_6^4(z,w),\,\, v=\psi^3_8(z,w)\}$, $T(B_{11})=D_4$;
\item[$12)$] $F_{12}=\{ u=\varphi_6(z,w),\,\, v=\psi_8(z,w)\}$, $T(B_{12})=E_6$;
\item[$13)$] $F_{13}=\{ u=\varphi^2_6(z,w),\,\, v=\theta_{12}(z,w)\}$, $T(B_{{13}})=E_7$;
\item[$14)$] $F_{14}=\{ u=\varphi_6(z,w),\,\, v=\theta_{12}^2(z,w)\}$, $T(B_{{14}})=A_7$;
\item[$15)$] $F_{15}=\{ u=\varphi^2_6(z,w),\,\, v=\theta_{12}^2(z,w)\}$, $T(B_{{15}})=D_6$;
\item[$16)$] $F_{16}=\{ u=\psi_{20}(z,w),\,\, v=\theta_{30}(z,w)\}$, $T(B_{{16}})=A_2$;
\item[$17)$] $F_{17}=\{ u=\psi_{20}(z,w),\,\, v=\theta_{30}^2(z,w)\}$, $T(B_{{17}})=A_5$;
\item[$18)$] $F_{18}=\{ u=\psi_{20}^3(z,w),\,\, v=\theta_{30}(z,w)\}$, $T(B_{{18}})=A_3$;
\item[$19)$] $F_{19}=\{ u=\varphi_{12}^5(z,w),\,\, v=\psi_{20}^3(z,w)\}$, $T(B_{19})=D_4$;
\item[$20)$] $F_{20}=\{ u=\varphi_{12}(z,w),\,\, v=\theta_{30}(z,w)\}$, $T(B_{{20}})=A_4$;
\item[$21)$] $F_{21}=\{ u=\varphi_{12}(z,w),\,\, v=\theta_{30}^2(z,w)\}$, $T(B_{{21}})=A_3$;
\item[$22)$] $F_{22}=\{ u=\varphi_{12}(z,w),\,\, v=\psi_{20}(z,w)\}$, $T(B_{22})=E_8$;
\end{itemize}

In case $1)$ the monodromy group of $F_1$ is $G_{F_1}=\mathbb S_3$; in cases $2_i)$, $i=1,2,3$, the groups $G_{F_{2_i,pq,q}}=G(pq,q,2)$ {\rm(}in notation used in \cite{Sh}{\rm )}; in cases $j)$, $3\leq j\leq 22$, the group $G_{F_j}$ is the group in \cite{Sh} with number $j$.
\end{thm}
\proof By Cartan's Lemma \cite{C}, the action of $G_{\widetilde F}$ on $\widetilde U$ can be linearized. Therefore $G_{\widetilde F}$ is a subgroup of $GL(2,\mathbb C)$ generated by reflections.
The list of germs of smooth Galois covers in Theorem \ref{Gal} is in one-to-one correspondence with the list in \cite{Sh} of finite subgroups of $GL(2,\mathbb C)$ generated by reflections. The singularity types of the germs of the branch curves $B_i$ of $F_i$, $i=1,\dots, 22$, are calculated in \cite{K}.
The rigidity of the covers $F_i$ easily follows from Cartan's Lemma. \qed
\subsection{Examples of finite groups which can (not) be realized as the monodromy groups of germs of finite covers.}
The following two statements are a direct consequence of Theorem \ref{Gal} and Proposition \ref{tras-subgr}.

\begin{cor} An abelian group $G$ can be realized as the monodromy group of the germ of a finite cover if and only if for each prime $p$, the number of cyclic $p$-primary factors entering in the presentation of $G$ as a direct product of cyclic groups is less than three.

The number of non-equivalent germs of finite covers
$F:(U,o')\to (V,o)$ whose monodrony group  $G_F=Ab$,
is equal to the number of non-isomorphic decompositions  $Ab=\mathbb Z_{n_1}\times \mathbb Z_{n_2}$ {\rm (}+ $1$, if $Ab$ is a cyclic group{\rm )}.\end{cor}
\proof There exists the unique relatively  prime subgroup of each abelian group $G$, namely, the trivial subgroup. \qed

\begin{cor} Groups $G=Q_8^k\times Ab$, $k\in \mathbb N$, where $$Q_8=\langle g_1, g_2 \mid g_1^2=g_2^2=[g_1,g_2],\,\, g_1^4=g_2^4=1\rangle$$ is the quaternion group and $Ab$ is any finite abelian group, are not realized as the monodromy groups of the germs of  finite covers.
\end{cor}
\proof There exists the unique relatively  prime subgroup of each group $G=Q_8^k\times Ab$, namely, the trivial subgroup and the groups  $G$ are not contained in the list of groups generated by reflections. \qed

\begin{prop} \label{alt4} The alternating group $\mathbb A_4$ can not be realized as the monodromy group of the germ of a finite cover.
\end{prop}
\proof By Proposition \ref{tras-subgr}, if $\mathbb A_4$ is the monodromy group of the germ $F:(U,o')\to (V,o)$ of some finite cover,
$\deg_{o'}F=d$, then $d$ is equal either $4$, or $6$, or $12$. By Corollary \ref{cor1}, $d\neq 4$ and by Theorem \ref{Gal}, $d\neq 12$, since the group $\mathbb A_4$ is not a group generated by reflections.

Assume that there exists the germ $F$ of a finite cover of $\deg_{o'}F=6$. Then the imbedding $G_F=\mathbb A_4\hookrightarrow \mathbb S_6$ corresponds to the action of $\mathbb A_4$ on the left cosets of a relatively simple subgroup $H_F\subset \mathbb A_4$ of order $2$. Consider the Galoisation
$\widetilde F:(\widetilde U,\widetilde o)\to (V,o)$ of the cover $F$.
The group $\mathbb A_4$ acts on $(\widetilde U,\widetilde o)$ such that $(\widetilde U,\widetilde o)/\mathbb A_4=(V,o)$ and
$(\widetilde U,\widetilde o)/H_1=(U,o')$, where $H_1=H_F$. The quotient cover
$F_{H_F}:(\widetilde U,\widetilde o)\to (U,o')$ is a two-sheeted germ of a Galois cover. Therefore it is branched along a germ of singular curve $B_{H_F}\subset U$, 
since $(\widetilde U ,\widetilde o)$ is a germ of singular (normal) surface by Theorem \ref{Gal}. Hence, the ramification curve $R_{H_F}\subset (\widetilde U,\widetilde o)$ of $\widetilde F$ is the germ of a singular curve also.

Note that the curve germ
$R_{H_F}$ can be defined as $$R_{H_F}=\{ p\in \widetilde U \mid g(p)=p\,\,\, \text{for}\,\, g\in H_F\}.$$

In the group $Kl_4\subset \mathbb A_4$ there exist three subgroups of order two, $H_1=H_F$, $H_2$, and $H_3$, conjugated in $\mathbb A_4$. Denote by $R_{i}=\{ p\in \widetilde U \mid g(p)=p\,\,\, \text{for}\,\, g\in H_i\}$, $i=2,3$. Note that $R_2\neq R_3$ and  $R_i$, $i=2,3$, are germs of singular curves.

The group $Kl_4\simeq H_1\times H_2$ acts on $(\widetilde U,\widetilde o)$. Consider the quotient germ $(\overline U,\overline o)=(\widetilde U ,\widetilde o)/Kl_4$ and the quotient mapping $\overline F:(\widetilde U ,\widetilde o)\to(\overline U,\overline o)$. We have $\overline F=F_{H_2}\circ F_{H_F}$, where $F_{H_2}: (U,o')\to (U,o')/H_2=(\overline U,\overline o)$, since $H_2$ is a normal subgroup of $Kl_4$ and therefore it acts on $(U,o')$. The mapping $F_{H_2}$ is ramified along the curve germ $F_{H_F}(R_2\cup R_3)$.
The germ  $F_{H_F}(R_2\cup R_3)\subset U$ is also a singular curve germ, since either $F_{H_F}(R_2)\neq F_{H_1}(R_3)$ or if $F_{H_F}(R_2)=F_{H_1}(R_3)$ then  $\deg F_{H_1\mid R_i}=1$,
since $\deg_{\widetilde o} F_{H_F}=2$.
On the other hand by Cartan's Lemma, the action of $H_2$ on $(U,o')$ can be linearized and the eigenvalues of the action of $H_2$ are  $\pm 1$. If two eigenvalues of the action of $H_2$ are $-1$, then the ramification locus of $F_{H_2}$ is the point $o'$, and if $H_2$ has the only one eigenvalue equals $-1$, when the ramification locus of $F_{H_2}$ is the germ of a non-singular curve.  A contradiction. \qed\\

Note that the groups $\mathbb S_2\simeq \mathbb Z_2$ and $\mathbb A_3\simeq \mathbb Z_3$ are the monodromy groups of cyclic Galois covers of degree $2$ and $3$.

\begin{prop} \label{alt2n-1} For $N\geq 3$ the symmetric group $\mathbb S_{N}$ and the alternating group $\mathbb A_{2N-1}$ are the monodromy groups of at least $[\frac{N-1}{2}]$ infinite series of rigid germs of finite covers.
\end{prop}
\proof Consider the germ $F_{m,n,k}:(U,o')\to (V,o)$ of a finite cover given by functions
\begin{equation} \label{last} \begin{array}{ll}
u=& z, \\
v= & \int{w^m(w-z^k)^{n}dw},\quad m,n,k\in \mathbb N.
\end{array}
\end{equation}
Obviously, $\deg_{o'} F_{m,n,k}=m+n+1$ and $F_{m,n,k}$ is ramified with multiplicity $m+1$ along the curve germ $R_1$ given by $w=0$ and with multiplicity $n+1$ along the curve germ $R_2$ given by $w-z^k=0$. Therefore the branch locus of $F_{m,n,k}$ is $B= B_1\cup B_2$, where $B_1=F_{m,n,k}(R_1)$ and $B_2=F_{m,n,k}(R_2)$.

We have
$$v= \displaystyle \sum_{i=0}^n\frac{(-1)^{n-i}{n\choose i}}{m+i+1}z^{k(n-i)} w^{m+i+1}.$$
It follows from (\ref{last}) that the restriction $F_{m,n,k |R_i}:R_i\to B_i$ of $F_{m,n,k}$ to $R_i$, $i=1,2$, is a biholomorphic mapping, $B_1$ is given by equation $v=0$ and $B_2$ is given by equation
$$v= \displaystyle \sum_{i=0}^n\frac{(-1)^{n-i}{n\choose i}}{m+i+1}u^{k(n+m+1)}.$$
Therefore the singularity type of the curve germ $B$ is $T(B)=A_{2k(n+m+1)-1}$.

By Lemma \ref{cl1.2}, the group $\pi_1^{loc}(B,o)$ is generated by two bypasses $\gamma_1$ around $B_1$ and $\gamma_2$ around $B_2$, and it follows from the above discussion that $F_*(\gamma_1)$ is a cycle of length $m+1$ and $F_*(\gamma_2)$ is a cycle of length $n+1$ in the symmetric group $\mathbb S_{m+n+1}$. Therefore, up to conjugation in $\mathbb S_{m+n+1}$, we can assume that $F_*(\gamma_1)=(1,2,\dots,m, m+n+1)$ and $F_*(\gamma_2)=(m+1,m+2,\dots, m+n+1)$, since $F_*(\gamma_1)$ and $F_*(\gamma_2)$ generate a transitive subgroup $G_{F_{m,n,k}}$ of $\mathbb S_{m+n+1}$.

\begin{lem}\label{alt} Let $G=\langle \tau,\sigma\rangle\subset \mathbb S_{m+n+1}$ be a subgroup generated by two cyclic permutations $\tau=(1,2,\dots,m, m+n+1)$ and $\sigma=(m+1,m+2,\dots, m+n+1)$, $m,n\in\mathbb N$. Then $\mathbb A_{m+n+1}\subset G$.
\end{lem}
\proof Without loss of generality, we can assume that $m\geq n$. Let us show that $G$ contains the all cycles of length 3.

We have $\varrho=\tau\sigma=(1,2,\dots, m+n+1)\in G$, 
\begin{equation} \label{q1} \eta_{m+1-i}=\tau^i\sigma\tau^{-i}= (m+1-i,m+1,\dots, m+n)\in G\,\, \text{for}\,\, i=1,\dots, m, \end{equation}
If $n=1$ then it follows from (\ref{q1}) that $G$ contains a transitive subgroup generated by transpositions $(j,m+1)$, $j=1,\dots, m$, and $\sigma$. Therefore $G=\mathbb S_{m+n+1}$. Therefore we can assume that $2\leq n\leq m$. If $m=2$, then direct check which we left to the reader shows that $\tau=(1,2,5)$ and $\sigma=(3,4,5)$ generate the group $\mathbb A_5$. Therefore we can assume that $m\geq 3$ and $m+n+1\geq 5$.

It is easy to check that
$$\begin{array}{rll} \eta_{j_1}\eta_{j_2}^{-1}= & (j_1,j_2,m+n)\in G & \quad  \text{for}\,\, 1\leq j_1<j_2\leq m, \\
(\eta_{j_1}\eta_{j_2}^{-1})(\eta_{j_1}\eta_{j_3}^{-1})^{-1}= & (j_1,j_2,j_3)\in G & \quad  \text{for}\,\, 1\leq j_1<j_2<j_3\leq m.\end{array}$$
Therefore
$$\varrho^{-k}[(\eta_{j_1}\eta_{j_2}^{-1})(\eta_{j_1}\eta_{j_3}^{-1})^{-1}]\varrho^{k}=  (j_1+k,j_2+k,j_3+k)\in G$$  for $1\leq j_1<j_2<j_3\leq m$ and $0\leq k\leq n+1$ and hence it suffices to prove that if the alternating groups $\mathbb A_k$ acting on the set $\{ 1,2, \dots, k\}$ and $\mathbb A_3$ acting on the set $\{ k-1,k,k+1\}$ are subgroups of $G\subset \mathbb S_{m+n+1}$, then the alternating group $\mathbb A_{k+1}$ acting on the set $\{ 1,2,\dots, k,k+1\}$ is also a subgroup of $G$. But, it is obvious, since, first, the cycles $(j_1,j_2,j_3)$ belong to $G$ for $\{ j_1,j_2,j_3\}\subset \{1,2,\dots, k\}$; second, the cycles $(j,k-1,k+1)$ and $(j,k,k+1)$ belong to $G$ for $1\leq j\leq k-2$, since the cycles
$(j,k-1,k)$ and $(j, k-1,k+1)$ belong to $G$; finally, the cycles $(j_1,j_2,k+1)$ belong to $G$ for $1\leq j_1<j_2\leq k-2$, since the cycles $(j_1,k,k+1)$ and $(j_2, k,k+1)$ belong to $G$. Therefore, to complete the proof of Lemma \ref{alt}, it suffices to notice that the subgroup $H_{k+1}\subset G$ generated by all cycles $(j_1,j_2,j_3)$ of length three,
$1\leq j_1<j_2<j_3\leq k+1$, is a normal subgroup of the group $\mathbb A_{k+1}$ if $k\geq 4$ and hence, $H_{k+1}=\mathbb A_{k+1}$. \qed \\

It follows from Lemma \ref{alt} that if one of the cycles $(1,2,\dots,m, m+n+1)$ and $(m+1,m+2,\dots, m+n+1)$ is an odd permutation, then these cycles generate the symmetric group $\mathbb S_{m+n+1}$, and if the both of these two cycles are even permutations, then they generate the group $\mathbb A_{m+n+1}$. In the first case, if we put $N=m+n+1$, then we obtain that there are at least $[\frac{N-1}{2}]$ infinite series ($k\in \mathbb N$) of the germs $F_{m,n,k}$ ($m$ or $n$ is an odd number) of finite covers with the monodromy group $G_{F_{m,n,k}}=\mathbb S_{N}$. In the second case, if we put $2N=m+n+2$, then we obtain that there are at least $\frac{N-1}{2}$ infinite series ($k\in \mathbb N$) of the germs $F_{m,n,k}$ ($m$ and $n$ are even numbers) of finite covers with the monodromy group
$G_{F_{m,n,k}}=\mathbb A_{2N-1}$.

The rigidity of the germs $F_{m,n,k}$ follows from the rigidity of the curve germs of singularity type $A_{2k(n+m+1)-1}$ and from Grauert - Remmert - Riemann - Stein Theorem, since the monodromy $F_{m,n,k*}$ is uniquely (up to conjugation in $\mathbb S_{m+n+1}$) determined by the triple $(m,n,k)$. \qed

\ifx\undefined\bysame
\newcommand{\bysame}{\leavevmode\hbox to3em{\hrulefill}\,}
\fi

\end{document}